\documentclass[12pt]{article}
\usepackage{amsmath}
\usepackage{amsfonts}

\advance \textheight by \topskip
\advance \textheight by 1pt
\makeatother
\def\bea{\begin{eqnarray}}
\def\ena{\end{eqnarray}}
\def\non{\nonumber}
\newcommand{\qed}{\hbox{\rule[-2pt]{3pt}{6pt}}}
\def\sgn{{\rm sgn}\,}
\def\bw{{\bar w}}

\newtheorem{prop}{Proposition}
\newtheorem{theorem}{Theorem}
\newtheorem{defn}{Definition}
\newtheorem{lemma}{Lemma}
\newtheorem{cor}{Corollary}

\def\pf{{\it Proof.}\,}
\def\wt{{\rm wt}\,}
\def\ord{{\rm ord}_\infty\,}
\def\hw{{\widehat \omega}}

\title{
Derivatives of Schur, Tau and Sigma Functions on Abel-Jacobi Images
}

\author{
Atsushi Nakayashiki\thanks{
e-mail: atsushi@tsuda.ac.jp}\\
Department of Mathematics,\\
Tsuda College\\
\\
Keijiro Yori\\
\\
\\
 \\
}

\date{}
\begin{document}
\maketitle
\centerline{
Dedicated to Michio Jimbo on his sixtieth birthday
}
\vskip15mm

\begin{abstract}
We study derivatives of Schur and tau functions from the view point
of the Abel-Jacobi map. We apply the results to establish several 
properties of derivatives of the sigma function of an $(n,s)$ curve.
As byproducts we have an expression of the prime form in terms of
derivatives of the sigma function and addition formulae which generalize
 those of Onishi for hyperelliptic sigma functions.
\end{abstract}
\clearpage

\section{Introduction}
The Riemann's theta function of an algebraic curve $X$ of genus $g$ 
can be considered, through the Abel-Jacobi map, as a multivalued
multiplicative analytic function on $X^g$. The Riemann's vanishing theorem
tells that the theta function shifted by the Riemann's constant vanishes 
identically on $X^{g-1}$. 
However it is possible to find certain derivatives of the theta
function such that they become multivalued multiplicative analytic functions 
on $X^{g-1}$. Onishi \cite{O1} found such derivatives explicitly
in the case of hyperelliptic curves. The extension of the results to the curve
$y^n=f(x)$ is given in \cite{MP}. 
These explicit derivatives of the theta function are used to 
construct certain addition formulae in \cite{O1}. 
The aim of this paper is to generalize and clarify the structure of 
the results on derivatives and
addition formulae in \cite{O1} by studying Schur and tau functions.

We consider a certain plane algebraic curve $X$, called
an $(n,s)$ curve \cite{BEL2}, which contains curves $y^n=f(x)$ 
as a special case. 
As in \cite{O1} we study sigma functions \cite{BEL1,N1} 
rather than Riemann's theta function since it is simpler 
to describe derivatives. 
Sigma functions can be expressed by the tau function of the KP-hierarchy 
\cite{EEG,EH,N2}. 
The expansion of the tau function with respect to Schur functions 
is known very explicitly due to Sato's theory of universal 
Grassmann manifold (UGM) \cite{SS,SN}. In the case corresponding to
the sigma function of an $(n,s)$ curve the expansion of the tau
function begins from the Schur function $s_{\lambda}(t)$ corresponding to
the partition $\lambda$ determined from the gap sequence at $\infty$ of $X$.
Notice that Schur functions themselves can be considered as a special
case of tau functions.

For a theta function solution of the KP-hierarchy the image of 
the Abel-Jacobi map of a point on a Riemann surface is transformed, in the tau function, to the vector of the form
\bea
&&
[z]=(z,z^2/2,z^3/3,...),
\label{eq-1-1}
\ena
where $z$ being a local coordinate at a base point.
Being motivated by this we consider, in general, the map 
$z\mapsto [z]$ as an analogue of the
Abel-Jacobi map for Schur and tau functions. For the Schur function
corresponding to an $(n,s)$ curve
 a similar map is considered in \cite{BEL2} as the rational limit of 
the Abel-Jacobi map.

The Schur function $s_\lambda(t)$, $t=(t_1,t_2,...)$, corresponding to a
partition $\lambda=(\lambda_1,...,\lambda_l)$ is the polynomial in 
$t_1,t_2,...$ defined by
\bea
&&
s_\lambda(t)=\det\left(p_{\lambda_i-i+j}(t)\right)_{1\leq i,j\leq l},
\quad
\exp({\sum_{i=1}^\infty t_ik^i})
=\sum_{i=1}^\infty p_i(t) k^i.
\non
\ena
We firstly study, for each $k$ satisfying $k\leq g$,
 the condition under which  a derivative
\bea
&&
\partial^\alpha s_\lambda([z_1]+\cdots+[z_k]),
\label{eq-1-2}
\ena
vanishes identically, where, for $\alpha=(\alpha_1,\alpha_2,...)$,
 $\partial^\alpha$ denote $\partial_1^{\alpha_1}\partial_2^{\alpha_2}\cdots$ and 
$\partial_i=\partial/\partial t_i$.
A sufficient condition can easily be found.
 Let us define the weight of $\alpha$
by $\wt\alpha=\sum_{i=1}^\infty i\alpha_i$ and set $N_{\lambda,k}=\lambda_{k+1}+\cdots+\lambda_l$. Then the derivative (\ref{eq-1-2}) vanishes, if $\wt\alpha<N_{\lambda,k}$. 

Concerning to derivatives such that (\ref{eq-1-2}) does not vanish identically
we have found two kinds of $\alpha$ satisfying $\wt\alpha=N_{\lambda,k}$. 
One is $\alpha=(N_{\lambda,k},0,0,...)$ for which the following recursive
relation holds:
\bea
&&
\partial_1^{N_{\lambda,k}}s_\lambda(\,\sum_{i=1}^k\,[z_i]\,)=
\frac{c'_{\lambda,k}}{c'_{\lambda,k-1}}
\partial_1^{N_{\lambda,k-1}}
s_\lambda(\,\sum_{i=1}^{k-1}\,[z_i]\,)z^{\lambda_k}+O(z_k^{\lambda_k+1}) ,
\label{eq-1-3}
\ena
where $c'_{\lambda,k}$ is a certain constant (Theorem \ref{th-2-4}).

The other kind of derivatives exist only for 
$\lambda$ corresponding to a gap sequence. 
A gap sequence of genus $g$ is a sequence of positive integers
$w_1<\dots<w_g$ such that its complement in the set of non-negative
integers ${\mathbb Z}_{\geq 0}$ is a semi-group. 
To each gap sequence a partition 
$\lambda=(\lambda_1,...,\lambda_g)$ is associated by 
\bea
&&
\lambda=(w_g,...,w_2,w_1)-(g-1,...,1,0).
\non
\ena
Let $w_1^\ast<w_2^\ast<\cdots$ be the complement of $\{w_i\}$
in ${\mathbb Z}_{\geq 0}$. For each $k$ the number $m_k$ and the
sequence $a^{(k)}_j$, $1\leq j\leq m_k$, are defined by
\bea
m_k&=&\sharp\{i|w_i^\ast<g-k\},
\non
\\
(a^{(k)}_1,...,a^{(k)}_{m_k})&=&(w_{g-k},w_{g-k-1},...,w_{g-k-m_k+1})-
(w_1^\ast,...,w_{m_k}^\ast).
\non
\ena
Then $\sum_{j=1}^{m_k}a^{(k)}_j=N_{\lambda,k}$ and the following
relation is valid:
\bea
&&
\partial_{a^{(k)}_1}\cdots \partial_{a^{(k)}_{m_k}}
s_\lambda(\,\sum_{i=1}^k\,[z_i]\,)=
\pm \partial_{a^{(k-1)}_1}\cdots \partial_{a^{(k-1)}_{m_k}}
s_\lambda(\,\sum_{i=1}^{k-1}\,[z_i]\,)z_k^{\lambda_k}+O(z_k^{\lambda_k+1}).
\label{eq-1-4}
\ena
 These derivatives generalize those of \cite{O1}\cite{MP}. Our construction
here clarifies the condition under which extensions of derivatives in \cite{O1} exist.

The tau function corresponding to a point of the cell $UGM^\lambda$ of UGM 
specified by a partition $\lambda$ has the expansion of the form
\bea
&&
\tau(t)=s_\lambda(t)+\sum_{\lambda<\mu} \xi_\mu s_\mu(t).
\label{eq-1-5}
\ena
We show that the vanishing property and the equations
(\ref{eq-1-3}), (\ref{eq-1-4}) for Schur functions hold without any
change if Schur functions are replaced by tau functions.
To this end we need to study 
derivatives of Schur functions $s_\mu(t)$ corresponding to partitions
$\mu$ satisfying $\lambda\leq \mu$ simultaneously. 
For example we have to study properties 
of "$a^{(k)}_j$-derivatives" of $s_\mu(t)$ where  $a^{(k)}_j$
are determined from $\lambda$. 

In the case corresponding to $(n,s)$ curves all the properties of tau functions established in this way are transplanted to sigma functions without much difficulty using the relation of the sigma function with
 the tau function. 

For applications to addition formulae we need to study derivatives
of Schur functions not only at $[z_1]+\cdots+[z_k]$ but at $[z_1]-[z_2]$.
In this case we have 
\bea
&&
\partial_1^{N'_{\lambda,1}}s_\lambda([z_1]-[z_2])=(-1)^{l-1}
\frac{c_\lambda}{c'_{\lambda,1}}\partial_1^{N_{\lambda,1}} 
s_\lambda([z_1])z_2^{l-1}+O(z_2^l),
\label{eq-1-6}
\ena
where $N'_{\lambda,1}=\lambda_2+\cdots+\lambda_l-l+1$ and 
 $c_\lambda$ is the constant given in Theorem \ref{th-2-2}.
It can be proved using the rational analogue of the 
Riemann's vanishing theorem for Schur functions \cite{BEL2}.  Again (\ref{eq-1-6}) and related properties are valid for
 tau and sigma functions without any change. 
As a corollary we obtain the expression of the prime form \cite{F,M2,N1} 
in terms of a certain derivative of the sigma function and consequently 
closed addition formulae for sigma functions. Here "closed" means
"without using prime form". 
The simplest example of the addition formula in the case of
an $(n,s)$ curve $X:y^n-\sum \lambda_{ij}x^iy^j=0$, is
\bea
&&
\frac{\partial_{u_1}^{N_{\lambda,2}}\sigma(p_2+p_1)
\partial_{u_1}^{N'_{\lambda,1}}\sigma(p_2-p_1)}
{(\partial_{u_1}^{N_{\lambda,1}}\sigma(p_1))^2
(\partial_{u_1}^{N_{\lambda,1}}\sigma(p_2))^2}
=
(-1)^g c_{\lambda}(c'_{\lambda,1})^{-4}c'_{\lambda,2}
(x_2-x_1),
\label{eq-1-7}
\ena
where $p_i\in X$ is identified with its Abel-Jacobi image,
$x_i=x(p_i)$ and $\lambda$ is the partition corresponding to
the gap sequence at $\infty$ of $X$. It generalizes the famous 
addition formula for Weierstrass' sigma function 
\bea
&&
\frac{\sigma(u_1+u_2)\sigma(u_1-u_2)}{\sigma(u_1)^2\sigma(u_2)^2}
=\wp(u_2)-\wp(u_1),
\non
\ena
since $(x_i,y_i)=(\wp(u_i),\wp'(u_i))$, $i=1,2$, are two points on $y^2=4x^3-g_2x-g_3$ and the right hand side can be written as $x_2-x_1$.
The formulae in \cite{O1} for hyperellitic sigma functions are recovered
if we use "$a^{(k)}_j$-derivatives" instead of $u_1$-derivative (see
the remark after Corollary \ref{cor-5-5}).

The present paper is organized as follows.
In section two properties of derivatives of Schur functions are studied.
The notion of gap sequence and the sequence $a^{(k)}_i$ are introduced. 
We lift the properties of Schur function in section two to functions
satisfying similar expansion to the tau functions of the KP-hierarchy
in section 3. In section 4 the properties on derivatives of the sigma 
function are proved using the sigma function expression of the tau function.
The expression of the prime form in terms of a derivative of the sigma function
of an $(n,s)$ curve is given in section 5. 
Addition formulae for sigma functions are proved.

\section{Schur function}
A sequence of non-negative integers $\lambda=(\lambda_1,...,\lambda_l)$
satisfying $\lambda_1\geq \cdots\geq \lambda_l$ is called a partition.
The number of non-zero elements in $\lambda$ is called the length of
$\lambda$ and is denoted by $l(\lambda)$. We identify $\lambda$ with
 partitions which are obtained from $\lambda$ by adding arbitrary
number of $0$'s, i.e. $(\lambda_1,...,\lambda_l,0,...,0)$. 
We set $|\lambda|=\lambda_1+\cdots+\lambda_l$.

Let $t=(t_1,t_2,t_3,...)$ and $p_n(t)$ the polynomial in $t$ defined
by
\bea
&&
\exp(\sum_{n=1}^\infty t_nk^n)=\sum_{n=0}^\infty p_n(t) k^n.
\label{pn}
\ena
We set $p_n(t)=0$ for $n<0$.

For a partition $\lambda=(\lambda_1,...,\lambda_l)$ Schur functions
$s_\lambda(t)$ and $S_\lambda(x)$ are defined by
\bea
s_\lambda(t)&=&\det(p_{\lambda_i-i+j}(t))_{1\leq i,j\leq l},
\non
\\
S_\lambda(x)&=&\frac{\det(x_j^{\lambda_i+l-i})_{1\leq i,j\leq l}}{
\prod_{i<j}(x_i-x_j)}.
\ena
The function $S_\lambda(x)$ is a symmetric polynomial of 
$x_1$,...,$x_l$ which is homogeneous of degree $|\lambda|$.

We introduce the symbol $[x]$ by
\bea
&&
[x]=(x,\frac{x^2}{2},\frac{x^3}{3},...),
\non
\ena
which is an analogue of Abel-Jacobi map in the theory of Schur functions.
With this symbol, $s_\lambda(t)$ and $S_\lambda(x)$ are related
by 
\bea
&&
s_\lambda(\sum_{i=1}^n[x_i])=S_\lambda(x),
\non
\ena
for $n\geq l(\lambda)$. From this relation we have

\begin{prop}\label{prop-2-1}
Let $\lambda=(\lambda_1,...,\lambda_l)$ be a partition of length $l$.
Then 
\vskip2mm
\noindent
(i) $\displaystyle{
s_\lambda\bigl(\,\sum_{i=1}^l\,[x_i]\,\bigr)=
s_{(\lambda_1,...,\lambda_{l-1})}\bigl(\,\sum_{i=1}^{l-1}\,[x_i]\,\bigr)x_l^{\lambda_l}
+O(x_l^{\lambda_l+1})}$.

\vskip2mm
\noindent
(ii) If $k<l$, $\displaystyle{
s_\lambda\bigl(\,\sum_{i=1}^k\,[x_i]\,\bigr)=0.}$
\end{prop}
\noindent
\pf
(i) It immediately follows from the definition of
$S_\lambda(x)$.

\noindent
(ii) We have
\bea
&&
s_\lambda\bigl(\,\sum_{i=1}^k\,[x_i]\,\bigr)=
s_\lambda\bigl(\,\sum_{i=1}^k\,[x_i]+[0]+\cdots+[0]\bigr).
\non
\ena
The right hand side is zero by (i) since $\lambda_l\geq 1$. \qed

Let $G^c$ be a subset of the set of non-negative integers 
${\mathbb Z}_{\geq 0}$. We assume that $G^c$ is a semi-group,
that is, it is closed under addition and contains $0$. Set
$G={\mathbb Z}_{\geq 0}\backslash G^c$.

\begin{defn}\label{def-gapseq}
Let $g$ be a positive integer.
 $G$ is called a gap sequence of genus $g$, if $\sharp G=g$.
Elements of $G$ and $G^c$ are called gaps and non-gaps respectively.
\end{defn}

For a gap sequence of genus $g$ enumerate elements of
$G$ and $G^c$ respectively  as
\bea
&&
w_1<w_2<\cdots<w_g,
\non
\\
&&
0=w_1^\ast<w_2^\ast<w_3^\ast<\cdots.
\non
\ena
Then $w_1=1$. For, otherwise $G^c$ contains $1$ and 
$G^c={\mathbb Z}_{\geq 0}$ which is impossible due to $g\geq 1$.
With this notation in mind we sometimes use $(w_1,...,w_g)$ to denote 
 a gap sequence instead of $\{w_1,...,w_g\}$.

\vskip2mm
\noindent
{\bf Example 1} Let $(n,s)$ be a pair of relatively prime 
integers such that $n,s\geq 2$. We set
\bea
&&
G^c=\{in+js\,|\, i,j\geq 0\}.
\non
\ena
Then $G$ is a gap sequence of genus $g=1/2(n-1)(s-1)$ \cite{BEL2}.
We call $G$ the gap sequence of type $(n,s)$. It is characterized 
by the condition that $G^c$ is generated by two elements.

\vskip2mm
\noindent
{\bf Example 2} Let $G=\{1,2,3,7\}$ and 
$G^c={\mathbb Z}_{\geq 0}\backslash G$. Then $G$ is a gap sequence
of genus four. In this case $G^c$ is generated by $4,5,6$. 
Therefore $G$ is not of type $(n,s)$ for any 
$(n,s)$.
\vskip2mm

In this way the gap sequences are classified by the minimum number
of generators of $G^c$.

For a gap sequence $\{w_1,...,w_g\}$ we associate
a partition $\lambda$ by 
\bea
&&
\lambda=(w_g,...,w_1)-(g-1,...,1,0).
\non
\ena

A special property of the partition determined from a 
gap sequence is the following.

\begin{prop}\label{prop-2-2}
If $\lambda$ is determined from a gap sequence
$(w_1,...,w_g)$, then $s_\lambda(t)$ does not depend on $t_i$, $i\notin\{w_1,...,w_g\}$.
\end{prop}
\vskip2mm
\noindent

In order to prove the proposition we introduce some notation.

For a partition $\lambda=(\lambda_1,...,\lambda_l)$ we associate
a strictly decreasing sequence of numbers $\bw_i$ by
\bea
&&
(\bw_1,...,\bw_l)=(\lambda_1,...,\lambda_l)+(l-1,l-2,...,0).
\non
\ena
By this correspondence the set of partitions of length at most $l$
bijectively corresponds to the set of strictly decreasing sequence
 non-negative integers $\bw_1>...>\bw_l\geq 0$. 

For $(\bw_1,...,\bw_l)$ we set
\bea
&&
(w_l,...,w_1)=(\bw_1,...,\bw_l).
\non
\ena
The introduction of the notation $\bw_i$ is for the sake of simplicity
in proofs and that of $w_i$ is for the sake of being consistent with
the notation of gap sequence.

For integers $i_1,...,i_l$ define the symbol $[i_1,...,i_l]$ as the determinant of the $l\times l$ matrix whose $j$-th row is
\bea
&&
\left(...,p_{i_j-1}(t),p_{i_j}(t)\right).
\non
\ena
We write $[i_1,...,i_l](t)$ if it is necessary to write explicitly 
 the dependence on $t$.

By the definition, $[i_1,...,i_l]$ is skew symmetric in the numbers
$i_1,...,i_l$ and becomes zero if two numbers coincide or some number 
is negative. 

With this notation 
\bea
&&
s_\lambda(t)=[\bw_1,...,\bw_l].
\non
\ena
Differentiating (\ref{pn}) by $t_i$ we have
\bea
&&
\partial_i p_n(t)=p_{n-i}(t),
\quad
\partial_i=\frac{\partial}{\partial t_i}.
\non
\ena
Therefore we have
\bea
&&
\partial_i s_\lambda(t)=\sum_{j=1}^l [\bw_1,...,\bw_j-i,...,\bw_l].
\non
\ena
\vskip3mm

\noindent
{\it Proof of Proposition \ref{prop-2-2}.}\par
We have to show, for $i\geq 1$,
\bea
&&
\partial_{w_i^\ast}s_\lambda(t)=\sum_{j=1}^lD_j=0,
\quad
D_j=[w_l,...,w_j-w_i^\ast,...,w_1].
\label{prop21-eq-1}
\ena
If $w_j-w_i^\ast<0$, obviously $D_j=0$. 
Suppose that $w_j-w_i^\ast>0$. Let $G=\{w_1,...,w_g\}$.
Then $w_j-w_i^\ast\in G$. For, if $w_j-w_i^\ast\in G^c$ then 
$w_j\in G^c+w_i^\ast\subset G^c$ which is absurd. Thus $w_j-w_i^\ast=w_k$
for some $k$. Notice that $w_i^\ast\geq 1$ and $k\neq j$, since $i\geq 1$.
Therefore $D_j=0$ because two rows coincide.  Consequently (\ref{prop21-eq-1})
is proved. \qed

\begin{defn}\label{a-sequence}
Let $G$ be a gap sequence of genus $g$. For $0\leq k\leq g-1$ we
define a positive integer $m_k$ and a sequence of 
integers $a^{(k)}_i$, $1\leq i\leq m_k$ by
\bea
&&
m_k=\sharp\{i\,|\, w_i^\ast<g-k\},
\non
\\
&&
(a^{(k)}_1,...,a^{(k)}_{m_k})=(w_{g-k},w_{g-k-1},...,w_{g-k-m_k+1})-
(w_1^\ast,...,w_{m_k}^\ast).
\non
\ena
\end{defn}
\vskip5mm

\noindent
{\bf Example 3} For the gap sequence of type $(2,2g+1)$ we have
\bea
&&
(w_1,w_2,...,w_g)=(1,3,...,2g-1),
\quad
(w_1^\ast,w_2^\ast,w_3^\ast,...)=(0,2,4,...).
\non
\ena
Then
\bea
&&
m_k=\sharp\{i\,|\,2i-2<g-k\}=\left[\frac{g-k+1}{2}\right],
\non
\\
&&
(a^{(k)}_1, a^{(k)}_2,...)=(2g-2k-1,2g-2k-5,2g-2k-9,...).
\non
\ena
This sequence recovers the rule for derivatives in \cite{O1}.
\vskip5mm

For a partition $\lambda=(\lambda_1,...,\lambda_l)$ 
 and a number $k$ such that $0\leq k\leq l-1$ we set
\bea
&&
N_{\lambda,k}=\lambda_{k+1}+\cdots+\lambda_l.
\label{new-eq-2-1}
\ena

\begin{lemma}\label{lem-2-1}
(i) $a^{(k)}_1>\cdots>a^{(k)}_{m_k}\geq 1$.
\vskip2mm
\noindent
(ii) Each $a^{(k)}_i$ belongs to $G$.
\vskip2mm
\noindent
(iii) Let $\lambda$ be the partition determined from $G$ then
\bea
&&
\sum_{i=1}^{m_k} a^{(k)}_i=N_{\lambda,k}.
\non
\ena
\end{lemma}
\vskip2mm
\noindent
\pf
(i) Notice that $(w_{g-k},w_{g-k-1},...)$ is strictly decreasing and $(w_1^\ast,w_2^\ast,...)$ is increasing. Therefore $\{a^{(k)}_i\}$ is
strictly decreasing. 
Since $G$ and $G^c$ are complement to each other we have
\bea
&&
\{0,1,...,g-k-1\}=\{w_1^\ast,...,w_{m_k}^\ast\}
\sqcup
\{w_1,...,w_{g-k-m_k}\}.
\label{eq-2-1}
\ena
Then, by the definition of the number $m_k$,
\bea
&&
w_1^\ast<\cdots<w_{m_k}^\ast<g-k \leq w_{{m_k}+1}^\ast<\cdots,
\non
\\
&&
w_1<\cdots<w_{g-k-m_k}<g-k\leq w_{g-k-m_k+1}<\cdots<w_{g-k}<\cdots.
\label{eq-2-2}
\ena
In particular $a^{(k)}_{m_k}=w_{g-k-m_k+1}-w^\ast_{m_k}\geq 1$.

\noindent
(ii) Suppose that $a^{(k)}_j\in G^c$. Since $G^c$ is a semi-group we have
\bea
&&
w_{g-k-j+1}=a^{(k)}_j+w_j^\ast\in G^c,
\non
\ena
which is absurd. Thus $a^{(k)}_j\in G$.

\noindent
(iii) By (\ref{eq-2-1}) we have
\bea
\sum_{i=1}^{m_k} a^{(k)}_i&=&
\sum_{i=g-k-m_k+1}^{g-k} w_i-\sum_{i=1}^{m_k} w_i^\ast
\non
\\
&=&
\sum_{i=g-k-m_k+1}^{g-k} w_i-\left(\sum_{i=1}^{g-k-1} i-\sum_{i=1}^{g-k-m_k}w_i\right)
\non
\\
&=&\sum_{i=1}^{g-k}w_i-\sum_{i=1}^{g-k-1}i=\sum_{i=k+1}^{g} \lambda_i.
\non
\ena
\qed
\vskip2mm

For $\alpha=(\alpha_1,\alpha_2,...)$ with finite number of non-zero
components we define the weight of $\alpha$ and the symbol $\partial^\alpha$ by
\bea
&&
\wt \alpha=\sum_{i=1}^\infty i \alpha_i,
\quad
\partial^\alpha=\partial_1^{\alpha_1}\partial_2^{\alpha_2}\cdots.
\non
\ena
The weight of $\partial^\alpha$ is defined to be the weight of $\alpha$.

\begin{prop}\label{prop-2-3}
 Let $\lambda=(\lambda_1,...,\lambda_l)$ be a partition and
$0\leq k\leq l-1$. If $\wt \alpha<N_{\lambda,k}$ we have
\bea
&&
\partial^\alpha s_\lambda\bigl(\,\sum_{i=1}^k\,[x_i]\,\bigr)=0.
\non
\ena
For $k=0$ the right hand side should be understood as 
$\partial^\alpha s_\lambda(0)$.
\end{prop}
\vskip2mm
\noindent
\pf
Notice that $\partial^\alpha s_\lambda(t)$ is a linear combination
of determinants of the form
\bea
&&
[\bw_1-r_1,...,\bw_l-r_l],
\quad
r_1+\cdots+r_l=\wt \alpha.
\label{eq-2-3}
\ena
If (\ref{eq-2-3}) is not zero, $\bw_i-r_i$ are all non-negative
and different. Thus there exists a permutation $(i_1,...,i_l)$
of $(1,...,l)$ such that
\bea
&&
\bw_{i_1}-r_{i_1}>\cdots>\bw_{i_l}-r_{i_l}\geq 0.
\non
\ena
Let $\mu$ be the partition corresponding to this strictly decreasing
sequence. Then 
\bea
&&
s_\mu(t)=[\bw_{i_1}-r_{i_1},...,\bw_{i_l}-r_{i_l}].
\non
\ena

If $l(\mu)>k$, $s_\mu\bigl(\,\sum_{i=1}^k\,[x_i]\,\bigr)=0$ by (ii) of Proposition \ref{prop-2-1}.

We prove that $l(\mu)\leq k$ is impossible if $\wt \alpha<N_{\lambda,k}$.
Suppose that $l(\mu)\leq k$. Then $\mu=(\mu_1,...,\mu_k,0,...,0)$ and
\bea
&&
\bw_{i_l}-r_{i_l}=0,
\quad
\bw_{i_{l-1}}-r_{i_{l-1}}=1,...
,\quad
\bw_{i_{k+1}}-r_{i_{k+1}}=l-k-1.
\non
\ena
Therefore
\bea
&&
r_{i_l}=\bw_{i_l},
\quad
r_{i_{l-1}}=\bw_{i_{l-1}}-1,
...,
\quad
r_{i_{k+1}}=\bw_{i_{k+1}}-(l-k-1),
\non
\ena
and we have
\bea
r_{i_l}+\cdots+r_{i_{k+1}}
&=&
\bw_{i_l}+\cdots+\bw_{i_{k+1}}-\left(1+2+\cdots+l-k-1\right)
\non
\\
&\geq&
\bw_l+\cdots+\bw_{k+1}-\left(1+2+\cdots+l-k-1\right).
\non
\ena
On the other hand 
\bea
r_{i_l}+\cdots+r_{i_{k+1}}
&\leq&
r_l+\cdots+r_1=\wt \alpha
\non
\\
&<&
\lambda_{k+1}+\cdots+\lambda_l
\non
\\
&=&
\bw_{k+1}+\cdots+\bw_l-\left(1+2+\cdots+l-k-1\right),
\non
\ena
which is a contradiction. Thus Proposition \ref{prop-2-3}
is proved. \qed
\vskip2mm

\begin{theorem}\label{th-2-1} 
Let $\lambda=(\lambda_1,...,\lambda_g)$ be the partition determined from a gap sequence of genus $g$,
$0\leq k\leq g$ and $a^{(k)}_j$ the associated sequence of numbers for $k\neq g$.We set $s_{(\lambda_1,...,\lambda_k)}
\bigl(\,\sum_{i=1}^k\,[x_i]\,\bigr)=1$ for
$k=0$ and $\partial_{a^{(k)}_1}\cdots\partial_{a^{(k)}_{m_k}}=1$ for 
$k=g$.
\vskip2mm
\noindent
(i) We have
\bea
&&
\partial_{a^{(k)}_1}\cdots\partial_{a^{(k)}_{m_k}}
s_\lambda\bigl(\,\sum_{i=1}^k\,[x_i]\,\bigr)=
c_k s_{(\lambda_1,...,\lambda_k)}\bigl(\,\sum_{i=1}^k\,[x_i]\,\bigr),
\non
\ena
where $c_k=\pm1$, $k\neq g$ is given by the sign of the permutation
\bea
&&
c_k=\sgn
\left(\begin{array}{cccccc}
w_1^\ast&...&w_{m_k}^\ast&w_{g-k-m_k}&...,&w_1\\
g-k-1&g-k-2&...&...&1&0\\
\end{array}\right),
\non
\ena
 and $c_g=1$.
\vskip2mm
\noindent
(ii) Let $\mu=(\mu_1,...,\mu_g)$ be a partition such that
$\mu_i=\lambda_i$ for $k+1\leq i\leq g$. Then
\bea
&&
\partial_{a^{(k)}_1}\cdots\partial_{a^{(k)}_{m_k}}
s_\mu\bigl(\,\sum_{i=1}^k\,[x_i]\,\bigr)=
c_k s_{(\mu_1,...,\mu_k)}\bigl(\,\sum_{i=1}^k\,[x_i]\,\bigr),
\non
\ena
where $c_k$ is the same as in (i).
\end{theorem}
\vskip5mm

\noindent
{\bf Remark 1} For the gap sequence of type $(n,s)$  it can be 
checked that the derivative determined from 
the sequence $a^{(k)}_j$ is the same as that found in \cite{MP}.
In that case (i) of Theorem \ref{th-2-1} is proved in that paper.
\vskip3mm

\begin{lemma}\label{lem-2-2}
Let $\lambda=(\lambda_1,...,\lambda_l)$ be a partition, $0\leq k\leq l-1$ and
$r_1,...,r_l$ non-negative integers.
Suppose that the following conditions:
\bea
&&
\sum_{i=1}^l r_i=N_{\lambda,k},
\label{eq-2-4}
\\
&&
[\bw_1-r_1,...,\bw_l-r_l]\bigl(\,\sum_{i=1}^k\,[x_i]\,\bigr)\neq 0.
\label{eq-2-5}
\ena
Then
\vskip2mm
\noindent
(i) We have $r_i=0$ for $1\leq i\leq k$.
\vskip2mm
\noindent
(ii) The sequence $(\bw_{k+1}-r_{k+1},...,\bw_l-r_l)$ is a permutation of
$(l-k-1,...,1,0)$.
\vskip2mm
\noindent
(iii) We have 
\bea
&&
[\bw_1-r_1,...,\bw_l-r_l]\bigl(\,\sum_{i=1}^k\,[x_i]\,\bigr)=c\, s_{(\lambda_1,...,\lambda_k)}\bigl(\,\sum_{i=1}^k\,[x_i]\,\bigr),
\non
\ena
where $c=\pm1$.
\end{lemma}
\vskip2mm
\noindent
\pf
By the assumption (\ref{eq-2-5}) there exists a permutation $(i_1,...,i_l)$ of $(1,...,l)$ and a partition
$\mu=(\mu_1,...,\mu_l)$ such that
\bea
&&
\bw_{i_1}-r_{i_1}>\cdots>\bw_{i_l}-r_{i_l}\geq 0,
\non
\\
&&
s_\mu(t)=[\bw_{i_1}-r_{i_1},\cdots,\bw_{i_l}-r_{i_l}],
\non
\ena
and $l(\mu)\leq k$ as in the proof of Proposition \ref{prop-2-3}.
In particular $\mu_i=0$ for $i\geq k+1$ which means
\bea
&&
\bw_{i_l}-r_{i_l}=0,\quad ...,\quad \bw_{i_{k+1}}-r_{i_{k+1}}=l-k-1.
\non
\ena
By a similar calculation to that in the proof of Proposition \ref{prop-2-3} we
have 
\bea
r_{i_{k+1}}+\cdots+r_{i_l}&=&\bw_{i_{k+1}}+\cdots+\bw_{i_l}-\left(1+2+\cdots+l-k-1\right)
\non
\\
&\geq& \bw_{k+1}+\cdots+\bw_l-\left(1+2+\cdots+l-k-1\right),
\label{eq-2-6}
\ena
and
\bea
r_{i_{k+1}}+\cdots+r_{i_l}&\leq& r_1+\cdots+r_l
\non
\\
&=& \lambda_{k+1}+\cdots+\lambda_l
\non
\\
&=& 
\bw_{k+1}+\cdots+\bw_{l}-\left(1+2+\cdots+l-k-1\right),
\label{eq-2-7}
\ena
where we use (\ref{eq-2-4}).
Therefore every inequalities in (\ref{eq-2-6}) and (\ref{eq-2-7}) are 
equalities. Then $r_{i_1}=\cdots=r_{i_k}=0$ by (\ref{eq-2-7}) and 
$(i_{k+1},...,i_l)$ is a permutation of $(k+1,...,l)$ by (\ref{eq-2-6}).
It, then, implies that $(i_1,...,i_k)$ is a permutation of $(1,...,k)$.

Since
\bea
&&
(\bw_{i_1}-r_{i_1},...,\bw_{i_l}-r_{i_l})=
(\bw_{i_1},...,\bw_{i_k},\bw_{i_{k+1}}-r_{i_{k+1}},...,\bw_{i_l}-r_{i_l})
\non
\ena
and it is strictly decreasing, $(i_1,...,i_k)=(1,...,k)$.
Thus
\bea
[\bw_{i_1}-r_{i_1},...,\bw_{i_l}-r_{i_l}]&=&
[\bw_{1},...,\bw_{k},l-k-1,...,1,0].
\label{eq-2-8}
\ena

\begin{lemma}
For a positive integer $m$ and a set of integers $i_1,...,i_k$ we have
\bea
&&
[i_1,...,i_k,m-1,...,1,0]=[i_1-m,...,i_k-m].
\non
\ena
\end{lemma}
\vskip2mm
\noindent
\pf
Expand the determinant at $m+k$-th row, $m+k-1$-th row,...,until $k+1$-st row
 successively and get the result. \qed

Applying the lemma to (\ref{eq-2-8}) we have
\bea
[\bw_{i_1}-r_{i_1},...,\bw_{i_l}-r_{i_l}]&=&[\bw_1-(l-k),...,\bw_k-(l-k)]
\non
\\
&=&
s_{(\lambda_1,...,\lambda_k)}(t).
\non
\ena
Since $(i_1,...,i_l)$ is a permutation of $(1,...,l)$,
\bea
&&
[\bw_{1}-r_{1},...,\bw_{l}-r_{l}]=\pm s_{(\lambda_1,...,\lambda_k)}(t).
\non
\ena
\qed
\vskip5mm

\noindent
{\it Proof of Theorem \ref{th-2-1}.}\par

In this proof we fix $k$ and denote $a^{(k)}_j$ simply by $a_j$.
Recall that
\bea
&&
s_\lambda(t)=[w_g,...,w_1].
\non
\ena
We compute the value of $\partial_{a_1}\cdots\partial_{a_{m_k}}s_\lambda(t)$ at
$t=t^{(k)}:=[x_1]+\cdots+[x_k]$.
\vskip2mm

\noindent
Step 1. We first consider the term for which the row labeled by 
$w_{g-k-(i-1)}$ is differentiated 
by $\partial_{a_i}$ for $1\leq i\leq m_k$. It is of the form
\bea
&&
A:=[w_g,..,w_{g-k+1},w_{g-k}-a_1,...,w_{g-k-(m_k-1)}-a_{m_k},w_{g-k-m_k},...,w_1].
\non
\ena
By the definition of $a_i$
\bea
&&
w_{g-k-(i-1)}-a_i=w_i^\ast.
\non
\ena
Therefore 
\bea
&&
A=[w_g,..,w_{g-k+1},w_1^\ast,...,w_{m_k}^\ast,w_{g-k-m_k},...,w_1].
\non
\ena
Using (\ref{eq-2-1}) we have
\bea
A&=&c_k[w_g,..,w_{g-k+1},g-k-1,...,1,0]
\non
\\
&=&
c_k s_{(\lambda_1,...,\lambda_k)}(t).
\non
\ena
\vskip2mm

\noindent
Step 2. We prove that the terms differentiated in a different way from that 
in Step 1 are zero at $t=t^{(k)}$.

By Lemma \ref{lem-2-1} (iii) and Lemma \ref{lem-2-2} (i) the term is zero at 
$t^{(k)}$ if some row corresponding to $w_i$, $g-k+1\leq i\leq g$, is 
differentiated. Therefore, for non-zero terms, only the last $g-k$ rows are
differentiated. 

So let us consider a term for which only some of last $g-k$ rows are differentiated. Notice that a term is zero if some row is differentiated more than once.
In fact some row corresponding to $w_j$ with $g-k-m_k+1\leq j\leq g-k$ is not differentiated
in this case. By (\ref{eq-2-2}) $w_j\geq g-k$. Consequently it is impossible
for the sequence $(w_{g-k},...,w_1)$ to be a permutation of $(g-k-1,...,1,0)$.
Then this term is zero at $t^{(k)}$ by Lemma \ref{lem-2-2} (ii).

As a consequence of the above argument we know that a term is zero if some row labeled by $w_j$ with $g-k-m_k+1\leq j\leq g-k$ is not differentiated. So let us consider
a term for which each row corresponding to $w_j$ with $g-k-m_k+1\leq j\leq g-k$
is differentiated exactly once.
We assume that the row corresponding $w_{g-k-(i-1)}$ is differentiated by
$\partial_{a_i}$ for $1\leq i<j$ with some $j\leq m_k$ and $\partial_{a_j}$
differentiates the row corresponding to $w_{g-k-(j'-1)}$ for some $j'$ 
with $j<j'$. We have 
\bea
&&
w_{g-k}-a_1=w_1^\ast,\quad...,\quad w_{g-k-(j-2)}-a_{j-1}=w_{j-1}^\ast,
\non
\ena
and
\bea
w_{g-k-(j'-1)}-a_j&=&w_{g-k-(j'-1)}-\left( w_{g-k-(j-1)}-w_j^\ast \right)
\non
\\
&=&
w_j^\ast-(w_{g-k-(j-1)}-w_{g-k-(j'-1)})<w_j^\ast.
\label{eq-2-9}
\ena

If $w_{g-k-(j'-1)}-a_j$ belongs to $G^c$, we have
\bea
&&
w_{g-k-(j'-1)}-a_j\in \{w_1^\ast,...,w_{j-1}^\ast\},
\non
\ena
by (\ref{eq-2-9}). Thus the term is zero since two rows coincide.

Suppose that $w_{g-k-(j'-1)}-a_j$ belongs to $G$.
Then
\bea
&&
w_{g-k-(j'-1)}-a_j\in \{w_1,...,w_{g-k-m_k}\},
\non
\ena
since $w_j^\ast<g-k$ and (\ref{eq-2-2}). In this case the term in consideration
is zero since again two rows coincide. 
Thus (i) of Theorem \ref{th-2-1} is proved.
\vskip2mm

\noindent
Step 3. We prove (ii) of Theorem \ref{th-2-1}.
Let $w_g'>\cdots>w_1'$ be the strictly decreasing sequence corresponding to 
$\mu$, that is,
\bea
&&
(w_g',...,w_1')=(\mu_1,...,\mu_g)+(g-1,...,1,0).
\non
\ena
By assumption $w_i=w'_i$ for $1\leq i\leq g-k$. Define ${w_i'}^\ast$, $i\geq 0$
by 
\bea
&&
\{{w_i'}^\ast\,|\,i\geq 0\}={\mathbb Z}_{\geq 0}\backslash\{w_i'\},
\non
\\
&&
0={w_1'}^\ast<{w_2'}^\ast<\cdots. 
\non
\ena
Then $w_i^\ast={w_i'}^\ast$ for $1\leq i\leq m_k$, since
\bea
\{w_1^\ast,...,w_{m_k}^\ast\}
\sqcup
\{w'_1,...,w'_{g-k-m_k}\}
&=&
\{w_1^\ast,...,w_{m_k}^\ast\}
\sqcup
\{w_1,...,w_{g-k-m_k}\}
\non
\\
&=&\{0,1,...,g-k-1\}.
\non
\ena
As a consequence the arguments in step 1 and step 2 are valid without 
any change if $w_i$, $w_i^\ast$ are replaced by $w'_i$, ${w'_i}^\ast$ respectively. \qed

Next we study properties of Schur functions with respect to $t_1$ derivative.

\begin{theorem}\label{th-2-2}
Let $\lambda=(\lambda_1,...,\lambda_l)$ be a partition, $(w_l,...,w_1)$ the corresponding strictly decreasing sequence and $0\leq k\leq l$. 
 Then
\bea
&&
\partial_1^{N_{\lambda,k}}
s_\lambda\bigl(\,\sum_{i=1}^k\,[x_i]\,\bigr)=
c'_{\lambda,k}s_{(\lambda_1,...,\lambda_k)}
\bigl(\,\sum_{i=1}^k\,[x_i]\,\bigr),
\non
\ena
where
\bea
&&
c'_{\lambda,k}=
\frac{N_{\lambda,k}!}{\prod_{i=1}^{l-k}w_i!}
\prod_{i<j}^{l-k}(w_j-w_i).
\non
\ena
\end{theorem}
\vskip2mm
\noindent
\pf We have
\bea
&&
s_\lambda(t)=[w_l,...,w_1].
\non
\ena
By Leibniz's rule
\bea
&&
\partial_1^{N_{\lambda,k}}s_\lambda(t)=
\sum_{r_1+\cdots+r_l=N_{\lambda,k}}
\frac{N_{\lambda,k}!}{r_1!\cdots r_l!}
[w_l-r_l,...,w_1-r_1].
\label{eq-2-13}
\ena
By Lemma \ref{lem-2-2}, if $[w_l-r_l,...,w_1-r_1](t^{(k)})\neq 0$ then $r_i=0$ 
for $l-k+1\leq i\leq l$, $(w_{l-k},...,w_1-r_1)$ is a permutation of 
$(l-k-1,...,1,0)$ and
\bea
&&
[w_l-r_l,...,w_1-r_1](t^{(k)})=\sgn\left(
\begin{array}{cccc}
w_{l-k}&\cdots&\cdots&w_1-r_1\\
l-k-1&\cdots&1&0\\
\end{array}
\right).
\non
\ena
In this case we can write
\bea
&&
w_i-r_i=\sigma(i-1),
\quad
1\leq i\leq l-k,
\non
\ena
for some $\sigma$ of an element of the symmetric group $S_{l-k}$ 
acting on $\{0,1,...,l-k-1\}$.
We define $1/n!=0$ for $n<0$ for the sake of convenience. Then
\bea
&&
\partial_1^{N_{\lambda,k}}s_\lambda(t^{(k)})=
A_{\lambda,k}s_{(\lambda_1,...,\lambda_k)}(t^{(k)}),
\non
\ena
where
\bea
&&
A_{\lambda,k}=\sum_{\sigma\in S_{l-k}}\sgn\sigma\,
\frac{N_{\lambda,k}!}{(w_1-\sigma(0))!\cdots(w_{l-k}-\sigma(l-k-1))!}.
\non
\ena
We have
\bea
\frac{A_{\lambda,k}}{N_{\lambda,k}!}&=&
\det\left(\frac{1}{(w_i-(j-1))!}\right)_{1\leq i,j\leq l-k}
\label{eq-2-10}
\\
&=&
\prod_{i=1}^{l-k}\frac{1}{w_i!}
\det\left(\prod_{m=0}^{j-2}(w_i-m)\right)_{1\leq i,j\leq l-k},
\label{eq-2-11}
\ena
where we set $\displaystyle{\prod_{m=0}^{j-2}(w_i-m)=1}$ for $j=1$.
Notice that the rule $1/n!=0$ for $n<0$ is taken into account in rewriting
(\ref{eq-2-10}) to (\ref{eq-2-11}), since, if $w_i-(j-1)<0$ then
$\prod_{m=0}^{j-2}(w_i-m)=0$.

Let us set 
\bea
&&
D=\det\left(\prod_{m=0}^{j-2}(w_i-m)\right)_{1\leq i,j\leq l-k}.
\non
\ena
Expanding $\displaystyle{\prod_{m=0}^{j-2}(w_i-m)}$ in $w_i$ 
we easily have
\bea
&&
D=\det\left(w_i^{j-1}\right)_{1\leq i,j\leq l-k}=\prod_{i<j}^{l-k}(w_j-w_i),
\non
\ena
and consequently 
\bea
&&
\frac{A_{\lambda,k}}{N_{\lambda,k}!}=
\frac{\prod_{i<j}^{l-k}(w_j-w_i)}{\prod_{i=1}^{l-k} w_i!}.
\non
\ena
\qed

In order to study addition formulae of sigma functions we need to study
properties of Schur functions at $t=[x_1]-[x_2]$.

For a partition $\lambda=(\lambda_1,...,\lambda_l)$ let 
$\lambda'=(\lambda'_{1},...,\lambda'_{l'})$ be the conjugate of $\lambda$, i.e. $\lambda_i'=\sharp\{j\,|\, \lambda_j\geq i\}$.

\begin{theorem}\label{th-2-3}
Let $\lambda=(\lambda_1,...,\lambda_l)$ be a partition of length $l$,
$\lambda'=(\lambda'_{1},...,\lambda'_{l'})$ and $\tilde{\lambda'}=(\lambda'_1-1,...,\lambda'_{l'}-1)$. Then
\bea
&&
s_\lambda([x]-\sum_{i=1}^{l'}[x_i])=(-1)^{N_{\lambda,1}}
s_{\tilde{\lambda'}}(\sum_{i=1}^{l'}[x_i])\prod_{j=1}^{l'}(x-x_j).
\non
\ena
\end{theorem}
\vskip2mm
\noindent
\pf
This theorem is essentially proved in the proof of Theorem 5.5 in \cite{BEL2}.
In \cite{BEL2} $\lambda$ is assumed to be the partition corresponding to the gap sequence of type $(n,s)$. In that case $\lambda=\lambda'$ and the assertion 
in this theorem is not stated. 
Here we give a proof since it is a key theorem for applications
to addition formulae.
For the notational simplicity we prove the assertion 
by interchanging $\lambda$ and $\lambda'$. All facts and notation concerning Schur and symmetric functions used in this proof can be found 
in \cite{Mac}

Let $e_i=e_i(x_1,...,x_m)$ be the elementary symmetric function:
\bea
&&
\prod_{i=1}^m(t+x_i)=\sum_{i=0}^m e_i t^{m-i}.
\label{th33-eq-0}
\ena
They satisfy the relation
\bea
&&
e_i(x_1,...,x_m)=e_i(x_1,...,x_{m-1})+x_m e_{i-1}(x_1,...,x_{m-1}).
\label{th33-eq-1}
\ena
In general, for a partition $\mu=(\mu_1,...,\mu_m)$, 
the following equation holds:
\bea
&&
S_{\mu'}(x_1,...,x_m)=\det(e_{\mu_i-i+j})_{1\leq i,j\leq m}.
\label{th33-eq-2}
\ena
Let ${\mathbf a}_j$ be the column vector defined by
\bea
&&
{\mathbf a}_j={}^t(e_{\lambda_1-1+j},e_{\lambda_2-2+j},...,e_{\lambda_l-l+j}),
\non
\ena
where $e_r=e_r(x,x_1,...,x_l)$.

By (\ref{th33-eq-1}), (\ref{th33-eq-2}) we have
\bea
s_{\lambda'}([x]+[x_1]+\cdots+[x_l])&=&S_{\lambda'}(x,x_1,...,x_l)
\non
\\
&=&
\det(e_{\lambda_i-i+j})_{1\leq i,j\leq l}
\non
\\
&=&
\det\left({\mathbf a}_1+x{\mathbf a}_0,{\mathbf a}_2+x{\mathbf a}_1,...,
{\mathbf a}_{l-1}+x{\mathbf a}_l\right)
\non
\\
&=&
\sum_{j=0}^l x^j \det\left({\mathbf a}_0,...,{\mathbf a}_{j-1},
{\mathbf a}_{j+1},...,{\mathbf a}_{l}\right)
\non
\\
&=&
\det\left(
\begin{array}{cccc}
1&-x&\cdots&(-x)^l\\
{\mathbf a}_0&{\mathbf a}_1&\cdots&{\mathbf a}_l\\
\end{array}
\right).
\label{th33-eq-3}
\ena

Let $p_r=\sum_{i=1}^l x_i^k$ be the power sum symmetric function, 
$\omega$, ${\widehat \omega}$ and $\iota$ the automorphisms of the ring
of symmetric polynomials in $x_1,...,x_l$ defined by
\bea
&&
{\widehat \omega}(p_r)=(-1)^r p_r,
\quad
\iota(p_r)=-p_r,
\quad
\omega=\iota\circ {\widehat \omega}.
\label{th33-eq-4}
\ena
Notice that ${\widehat \omega}$ is, in terms of $x_j$, the map sending
$x_j$ to $-x_j$ for $1\leq j\leq l$.
Then
\bea
&&
s_{\lambda'}\bigl([x]-\sum_{i=1}^l\,[x_i]\bigr)=
(-1)^{|\lambda|}\omega\Bigl(
s_{\lambda'}\bigl([-x]+\sum_{i=1}^l\,[x_i]\bigr)
\Bigr).
\label{th33-eq-5}
\ena
It can be checked by computing the right hand side using 
(\ref{th33-eq-4}) and the relation 
$S_{\mu}(-x_1,...,-x_m)=(-1)^{|\mu|}S_{\mu}(x_1,...,x_m)$.

Let $h_i=h_i(x_1,...,x_l)$ be the complete symmetric function:
\bea
&&
\frac{1}{\prod_{i=1}^l(1-tx_i)}=
\sum_{i=0}^\infty h_i x^i.
\non
\ena
Then $\omega(e_i)=h_i$ and 
\bea
&&
\omega({\mathbf a}_j)={}^t(h_{\lambda_1-1+j},...,h_{\lambda_l-l+j}).
\label{th33-eq-6}
\ena
By (\ref{th33-eq-3}) and (\ref{th33-eq-5}) we have
\bea
&&
s_{\lambda'}([x]-[x_1]-\cdots-[x_l])=
(-1)^{|\lambda|}
\det\left(
\begin{array}{cccc}
1&x&\cdots&x^l\\
\omega({\mathbf a}_0)&\omega({\mathbf a}_1)&\cdots&\omega({\mathbf a}_l)\\
\end{array}
\right).
\label{th33-eq-7}
\ena
Using the relation,
\bea
&&
\sum_{j=0}^k (-1)^j e_j h_{k-j}=0,\quad k\geq 1,
\ena
we have
\bea
&&
\sum_{j=0}^k (-1)^j e_j \omega({\mathbf a}_{l-j})={\mathbf o}.
\label{th33-eq-8}
\ena
By (\ref{th33-eq-0}), (\ref{th33-eq-6}), (\ref{th33-eq-7}), (\ref{th33-eq-8})
we obtain
\bea
s_{\lambda'}\bigl([x]-\sum_{i=1}^l\,[x_i]\bigr)&=&
(-1)^{l+|\lambda|}
\det\left(
\omega({\mathbf a}_0),...,\omega({\mathbf a}_{l-1})\right)
\prod_{j=1}^l(x-x_j)
\non
\\
&=&
(-1)^{N_{\lambda',1}}
\det(h_{\lambda_i-1-i+j})_{1\leq i,j\leq l}
\prod_{j=1}^l(x-x_j).
\non
\ena
Then the theorem follows from
\bea
&&
S_{(\mu_1,...,\mu_m)}(x_1,...,x_m)=
\det(h_{\mu_i-i+j})_{1\leq i,j\leq m}.
\non
\ena
\qed

\begin{cor}\label{cor-2-1}
Let $\lambda=(\lambda_1,...,\lambda_l)$ be a partition of length $l$.Then
$s_\lambda([x_1]-[x_2])$ is not identically zero if and only if $\lambda_i=1$
for $2\leq i\leq l$, that is, $\lambda$ is a hook.
\end{cor}
\vskip2mm
\noindent
\pf Setting $x_i=0$ for $2\leq i\leq l'$ in Theorem \ref{th-2-3} we have
\bea
&&
s_\lambda([x]-[x_1])=(-1)^{N_{\lambda,1}}
s_{\tilde{\lambda'}}([x_1])x^{l'-1}(x-x_1).
\label{eq-2-12}
\ena
Thus $s_\lambda([x]-[x_1])\neq 0$ is equivalent to 
$s_{\tilde{\lambda'}}([x_1])\neq 0$. The latter is equivalent to
the condition that the length of $\tilde{\lambda'}$ is one.
It means that $\lambda'=(\lambda_1',1^{l'-1})$ which is equivalent to 
that $\lambda$ is a hook. \qed

\begin{theorem}\label{th-2-4}
Let $\lambda=(\lambda_1,...,\lambda_l)$ be a partition of length $l$,
$(w_l,...,w_1)$ the corresponding sequence and $N'_{\lambda,1}=\sum_{i=2}^l\lambda_i-l+1$. 
\vskip2mm
\noindent
(i) If $n<N'_{\lambda,1}$ 
\bea
&&
\partial_1^ns_\lambda([x_1]-[x_2])=0.
\non
\ena
\noindent
(ii) We have
\bea
&&
\partial_1^{N'_{\lambda,1}}s_\lambda([x_1]-[x_2])=
c_{\lambda}s_{(\lambda_1,1^{l-1})}([x_1]-[x_2]),
\non
\ena
where
\bea
&&
c_\lambda=\frac{N'_{\lambda,1}!}{\prod_{i=1}^{l-1}(w_i-1)!}
\prod_{i<j}^{l-1}(w_j-w_i).
\non
\ena
\noindent
(iii) Let $\mu=(\mu_1,...,\mu_{l'})$ be a partition of length $l'\geq l$ 
such that $\mu_i=\lambda_i$ for $2\leq i\leq l$ and $\mu_i=1$ for $i>l$.
Then
\bea
\partial_1^{N'_{\lambda,1}}s_\mu([x_1]-[x_2])=
c_{\lambda}s_{(\mu_1,1^{l'-1})}([x_1]-[x_2]).
\non
\ena
\noindent
(iv) For $m,n\geq 1$ we have
\bea
&&
a_{(m,1^{n-1})}([x_1]-[x_2]) =(-1)^{n-1}x_1^{m-1}x_2^{n-1}(x_1-x_2).
\non
\ena
\end{theorem}
\vskip2mm

\noindent
\pf Notice that
\bea
&&
\partial_1s_\lambda(t)=\sum_{i=1}^l[w_l,...,w_i-1,...,w_1].
\non
\ena
In the right hand side $[w_l,...,w_i-1,...,w_1]\neq 0$ if and only if all its components are different. In terms of the diagram of $\lambda$,  $\partial_1s_\lambda(t)$ is a sum of $s_\mu(t)$ with $\mu$ being the diagram obtained from
$\lambda$ by removing one box. For example
\bea
&&
\partial_1s_{(2,2,1)}(t)=s_{(2,1,1)}(t)+s_{(2,2)}(t).
\non
\ena
\vskip2mm
\noindent
(i) Notice that $N'_{\lambda,1}$ is a number of boxes on two to $l$-th rows 
 of the diagram of $\lambda$ which are on the right of the first column. 
Thus if $n<N'_{\lambda,1}$ it is impossible to get the
hook diagram by removing $n$ boxes from $\lambda$. Then the assertion of (i) follows from Corollary \ref{cor-2-1}.
\vskip2mm
\noindent
(ii) There is only one hook diagram in diagrams obtained from $\lambda$
by removing $N'_{\lambda,1}$ boxes. It is $\mu:=(\lambda_1,1^{l-1})$.
Let us compute the coefficient $c$ of $s_\mu(t)$ in 
$\partial_1^{N'_{\lambda,1}}s_\lambda(t)$.
Consider Equation (\ref{eq-2-13}) with $N_{\lambda,k}$ being replaced 
by $N'_{\lambda,1}$. In the right hand side, $s_\mu(t)$
appears only as a term such that
 $r_l=0$ and $(w_{l-1}-r_{l-1},...,w_1-r_1)$ is 
a permutation of $(l-1,...,2,1)$. Let us write, for $1\leq i\leq l-1$,
\bea
&&
w_i-r_i=\sigma(i), \quad \sigma\in S_{l-1}.
\non
\ena
Then by a similar calculation to that in the proof of Theorem \ref{th-2-2}
we have
\bea
&&
\frac{c}{N'_{\lambda,1}!}=
\sum_{\sigma\in S_{l-1}}\frac{\sgn\,\sigma}
{(w_1-\sigma(1))!\cdots(w_{l-1}-\sigma(l-1))!}
=
\frac{\prod_{i<j}^{l-1}(w_j-w_i)}{\prod_{i=1}^{l-1}(w_i-1)!}.
\non
\ena
\noindent
(iii) Similarly to the proof of (ii) the only Schur function appearing in
$\partial_1^{N'_{\lambda,1}}s_\mu(t)$ which does not vanish at $t=[x_1]-[x_2]$
is $s_\nu(t)$, $\nu=(\mu_1,1^{l'-1})$. Let us compute the coefficient $c'$
of $s_\nu(t)$ in $\partial_1^{N'_{\lambda,1}}s_\mu(t)$.

Let $(w'_{l'},...,w'_1)$ be the strictly decreasing sequence corresponding to
$\mu$. Then
\bea
&&
\partial_1^{N'_{\lambda,1}}s_\mu(t)=
\sum_{r_1+\cdots+r_{l'}=N'_{\lambda,1}}
\frac{N'_{\lambda,1}!}{r_1!\cdots r_{l'}!}
[w'_{l'}-r_l,...,w'_1-r_1].
\label{eq-2-14}
\ena
In the right hand side $[w'_{l'}-r_l,...,w'_1-r_1]$ is proportional to
$s_\nu(t)$ if and only if $r_i=0$ for $i=l'$ or $i<l'-l$, and 
$(w'_{l'-1}-r_{l'-1},...,w'_{l'-l+1}-r_{l'-l+1})$ is a permutation of
$(l'-1,l'-2,...,l'-l+1)$. 
Let us write
\bea
&&
w'_i-r_i=\sigma(i),
\quad
l'-l+1\leq i\leq l'-1,
\quad
\sigma\in S_{l-1}.
\non
\ena
Then
\bea
\frac{c'}{N'_{\lambda,1}!}&=&
\sum_{\sigma\in S_{l-1}}\frac{\sgn\,\sigma}{(w'_{l'-l+1}-\sigma(l'-l+1))!
\cdots (w'_{l'-1}-\sigma(l'-1))!}
\non
\\
&=&
\frac{\prod_{l'-l+1\leq i<j\leq l'-1}(w'_j-w'_i)}{\prod_{i=l'-l+1}^{l'-1}(w_i'-l'+l-1)!}.
\label{eq-2-15}
\ena

Let us rewrite $c'$ in terms of $\lambda_j$. 
By assumption $\mu_i=\lambda_i$ for $2\leq i\leq l$ which implies
\bea
&&
w_i'=\mu_{l'+1-i}+i-1=\lambda_{l'+1-i}+i-1,
\quad
l'-l+1\leq i\leq l'-1.
\non
\ena
Substitute it into (\ref{eq-2-15}) and get
\bea
&&
c'=\frac{N'_{\lambda,1}!}{\prod_{i=2}^l(\lambda_i+l-1-i)!}
\prod_{2\leq i<j\leq l}(\lambda_i-\lambda_j+j-i),
\non
\ena
which equals to $c_\lambda$.
\vskip2mm
\noindent
(iv) Set $\lambda=(m,1^{n-1})$ in (\ref{eq-2-12}). Then, using
$s_{(r)}([x])=x^r$, we get the assertion of (iv). \qed

\section{$\tau$-function}
In this section we lift the properties of Schur functions which have been
proved in the previous section to $\tau$-functions.

Let $\leq$ be the partial order on the set of partitions defined as
follows. For two partitions $\lambda=(\lambda_1,...,\lambda_l)$,
$\mu=(\mu_1,...,\mu_{l'})$, $\lambda\leq \mu$ if and only if
$\lambda_i\leq \mu_i$ for all $i$. 

For a partition $\lambda=(\lambda_1,...,\lambda_l)$ we consider
a function $\tau(t)$ given as a series of the form
\bea
&&
\tau(t)=s_\lambda(t)+\sum_{\lambda<\mu}\xi_\mu s_\mu(t),
\label{eq-3-1}
\ena
where $\xi_u\in {\mathbb C}$. 
\vskip2mm

\noindent
{\bf Example}. Let $X$ be a compact Riemann surface of genus $g\geq 1$,
$p_\infty$ a point of $X$, $w_1<\cdots<w_g$ the gap sequence at $p_\infty$ and $z$ a local coordinate at $p_\infty$. Embed the affine ring of $X\backslash \{p_\infty\}$ into Sato's  universal Grassmann manifold (UGM)
as in the paper \cite{N2}. Then the tau function corresponding to this 
point of UGM has the expansion of the form (\ref{eq-3-1}).

\begin{prop}\label{prop-3-1}
Let $\lambda=(\lambda_1,...,\lambda_l)$ a partition, $\tau(t)$ be a
function of the form (\ref{eq-3-1}) and $0\leq k\leq l-1$. Then,
if $\wt\alpha<N_{\lambda,k}$
\bea
&&
\partial^\alpha \tau\bigl(\,\sum_{i=1}^k\,[x_i]\,\bigr)=0.
\non
\ena
\end{prop}
\noindent
\pf For $\mu=(\mu_1,...,\mu_{l'})$ satisfying $\lambda\leq\mu$ we have
\bea
&&
\wt\alpha<\sum_{i=k+1}^l\lambda_i\leq \sum_{i=k+1}^{l'} \mu_i.
\non
\ena
Thus
\bea
&&
\partial^\alpha s_\mu\bigl(\,\sum_{i=1}^k\,[x_i]\,\bigr)=0,
\non
\ena
by Proposition \ref{prop-2-3}.
The assertion of the proposition follows from (\ref{eq-3-1}). \qed

For a function $\tau(t)$ of the form (\ref{eq-3-1}) and $1\leq k\leq l$
let $\tau^{(k)}(t)$ be the function defined by
\bea
&&
\tau^{(k)}(t)=s_{(\lambda_1,...,\lambda_k)}(t)+
\sum_{\mu}\xi_\mu s_{(\mu_1,...,\mu_k)}(t),
\non
\ena
where the sum in the right hand side is over all partitions 
$\mu=(\mu_1,...,\mu_l)$ such that $\lambda<\mu$ and 
$\mu_i=\lambda_i$ for $k+1\leq i\leq l$.
In particular $\tau^{(l)}(\sum_{i=1}^k[x_i])=\tau(\sum_{i=1}^k[x_i])$. 
We set $\tau^{(0)}(t)=1$.

\begin{theorem}\label{th-3-1}
Let $\lambda=(\lambda_1,...,\lambda_g)$ be the partition determined from a gap sequence 
of genus $g$, $\tau(t)$ a function of the form (\ref{eq-3-1}).
\vskip2mm
\noindent
(i) We have, for $0\leq k\leq g$, 
\bea
&&
\partial_{a^{(k)}_1}\cdots\partial_{a^{(k)}_{m_k}}
\tau\bigl(\,\sum_{i=1}^k\,[x_i]\,\bigr)=
c_k \tau^{(k)}\bigl(\,\sum_{i=1}^k\,[x_i]\,\bigr),
\non
\ena
where $c_k$ is the same as in Theorem \ref{th-2-1}.
\vskip2mm
\noindent
(ii) We have, for $k\geq 1$, 
\bea
&&
\tau^{(k)}\bigl(\,\sum_{i=1}^k\,[x_i]\,\bigr)=\tau^{(k-1)}\bigl(\,\sum_{i=1}^{k-1}\,[x_i]\,\bigr)x_k^{\lambda_k}+O(x_k^{\lambda_k+1}).
\non
\ena
\end{theorem}
\vskip2mm
\noindent
\pf
Let $\mu=(\mu_1,..,\mu_l)$ be a partition of length $l$ such that $\lambda\leq \mu$. Then $l\geq g$ and
\bea
&&
\wt\left(\partial_{a^{(k)}_1}\cdots\partial_{a^{(k)}_{m_k}}\right)
=\sum_{i=1}^{m_k}a^{(k)}_i
= \sum_{i=k+1}^{g} \lambda_i
\leq \sum_{i=k+1}^l \mu_i.
\label{eq-3-2}
\ena
If the inequality in the right hand side is not an equality, 
\bea
&&
\partial_{a^{(k)}_1}\cdots\partial_{a^{(k)}_{m_k}}
s_\mu\bigl(\,\sum_{i=1}^k\,[x_i]\,\bigr)=0,
\label{eq-3-3}
\ena
by Proposition \ref{prop-2-3}. Therefore, if the left hand side
of (\ref{eq-3-3}) does not vanish, $l=g$ and 
\bea
&&
\sum_{i=k+1}^g \mu_i=\sum_{i=k+1}^g \lambda_i.
\non
\ena
Since $\lambda_i\leq \mu_i$ for any $i$, it implies $\mu_i=\lambda_i$
for $k+1\leq i\leq g$. For such $\mu$ we have, by Theorem \ref{th-2-1},
\bea
&&
\partial_{a^{(k)}_1}\cdots\partial_{a^{(k)}_{m_k}}
s_\mu\bigl(\,\sum_{i=1}^k\,[x_i]\,\bigr)=
c_k s_{(\mu_1,...,\mu_k)}\bigl(\,\sum_{i=1}^k\,[x_i]\,\bigr).
\non
\ena
The assertion (i) follows from this.
\vskip2mm
\noindent
(ii) The assertion easily follows from (i) of Proposition \ref{prop-2-1}  and
the definition of $\tau^{(k)}(t)$. \qed

Combining (i) and (ii) of Theorem \ref{th-3-1} we have

\begin{cor}\label{cor-3-1} Under the same assumption as in Theorem 
\ref{th-3-1} we have, for $1\leq k\leq g$, 
\bea
&&
\partial_{a^{(k)}_1}\cdots\partial_{a^{(k)}_{m_k}}
\tau\bigl(\,\sum_{i=1}^k\,[x_i]\,\bigr)=
\frac{c_k}{c_{k-1}}
\partial_{a^{(k-1)}_1}\cdots\partial_{a^{(k-1)}_{m_{k-1}}}
\tau\bigl(\,\sum_{i=1}^{k-1}\,[x_i]\,\bigr)x_k^{\lambda_k}+
O(x_k^{\lambda_k+1}).
\non
\ena
\end{cor}

Corresponding Theorem \ref{th-2-2} we have

\begin{theorem}\label{th-3-2}
Let $\lambda=(\lambda_1,...,\lambda_l)$ be a partition
of length $l$, $\tau(t)$ a function of the form (\ref{eq-3-1}),
$0\leq k\leq l$. Then
\bea
&&
\partial_1^{N_{\lambda,k}}
\tau\bigl(\,\sum_{i=1}^k\,[x_i]\,\bigr)
=c'_{\lambda,k}\tau^{(k)}\bigl(\,\sum_{i=1}^{k}\,[x_i]\,\bigr),
\non
\ena
where $c'_{\lambda,k}$ is the same as in Theorem \ref{th-2-2}.
\end{theorem}
\vskip2mm
\noindent
\pf The theorem can be proved in a similar manner to Theorem \ref{th-3-1}
using Theorem \ref{th-2-2}. \qed

\begin{cor}\label{cor-3-2}
Under the same assumption as in Theorem \ref{th-3-2} we have,
for $1\leq k\leq l$, 
\bea
&&
\partial_1^{N_{\lambda,k}}\tau\bigl(\,\sum_{i=1}^k\,[x_i]\,\bigr)
=\frac{c'_{\lambda,k}}{c'_{\lambda,k-1}}
\partial_1^{N_{\lambda,k-1}}\tau\bigl(\,\sum_{i=1}^{k-1}\, [x_i]\,\bigr)x_k^{\lambda_k}+O(x_k^{\lambda_k+1}).
\non
\ena
\end{cor}

In order to state the properties for $\tau(t)$ corresponding to Theorem 
\ref{th-2-4} let us introduce one more function $\tau_2(t)$ associated 
with $\tau(t)$ by
\bea
&&
\tau_2(t)=s_{(\lambda_1,1^{l-1})}(t)+\sum_{\mu}\xi_\mu s_{(\mu_1,1^{l'-1})}(t),
\label{eq-3-4}
\ena
where the sum in the right hand side is over all partitions 
$\mu=(\mu_1,...,\mu_{l'})$ of length $l'\geq l$ satisfying $\lambda<\mu$, $\mu_i=\lambda_i$ for $2\leq i\leq l$ and $\mu_i=1$ for $i>l$.

\begin{theorem}\label{th-3-3}
Let $\lambda=(\lambda_1,...,\lambda_l)$ be a partition of length $l$ and 
$\tau(t)$ a function of the form (\ref{eq-3-1}).
\vskip2mm
\noindent
(i) If $n<N'_{\lambda,1}$
\bea
&&
\partial_1^n\tau([x_1]-[x_2])=0.
\non
\ena
\noindent
(ii) We have
\bea
&&
\partial_1^{N'_{\lambda,1}}\tau([x_1]-[x_2])=
c_{\lambda}\tau_2([x_1]-[x_2]),
\non
\ena
where $c_\lambda$ is the same as in Theorem \ref{th-2-4}.
\vskip2mm
\noindent
(iii) We have
\bea
&&
\tau_2([x_1]-[x_2])=
(-1)^{\lambda_1-1}x_1^{\lambda_1-1}x_2^{l-1}(x_1-x_2)(1+\cdots),
\non
\ena
where $\cdots$ part is a series in $x_1$, $x_2$ containing only 
terms proportional to $x_1^ix_2^j$ with $i+j>0$.
\vskip2mm
\noindent
(iv) We have the expansion
\bea
&&
\tau_2([x_1]-[x_2])=(-1)^{l-1}\tau^{(1)}([x_1])x_2^{l-1}+O(x_2^l).
\non
\ena
\end{theorem}
\vskip2mm
\noindent
\pf 
(i) By (i) of Theorem \ref{th-2-4} we have 
$\displaystyle{\partial_1^n s_\lambda([x_1]-[x_2])=0}$.

Suppose that $\lambda<\mu$ and $\mu=(\mu_1,...,\mu_{l'})$ is of length $l'$.
Then, $l'\geq l$ and 
\bea
&&
n<\sum_{i=2}^l \lambda_i-(l-1)\leq 
\sum_{i=2}^l \mu_i+\sum_{i=l+1}^{l'} (\mu_i-1)-(l-1)=
\sum_{i=2}^{l'}\mu_i-(l'-1).
\non
\ena
Thus $\displaystyle{\partial_1^n s_\mu([x_1]-[x_2])=0}$ by Theorem \ref{th-2-4} (i) and the assertion (i) is proved.
\vskip2mm
\noindent
(ii) By (ii) of Theorem \ref{th-2-4} we have
\bea
&&
\partial_1^{N'_{\lambda,1}}\tau([x_1]-[x_2])=
c_{\lambda}s_{(\lambda_1,1^{l-1})}([x_1]-[x_2])
+\sum_\mu\xi_\mu \partial_1^{N'_{\lambda,1}}s_\mu([x_1]-[x_2]).
\label{eq-3-5}
\ena
Let us compute the second term in the right hand side of (\ref{eq-3-5}).

Suppose that $\mu>\lambda$, $\mu=(\mu_1,...,\mu_{l'})$ is of 
length $l'$ and
$\partial_1^{N'_{\lambda,1}}s_\mu([x_1]-[x_2])\neq 0$. 
In such a case, similarly to the proof of Theorem \ref{th-2-4} (ii), 
it can be shown that $\mu$ should be of the form 
$\mu=(\mu_1,\lambda_2,...,\lambda_l,1^{l'-l})$.
Then
\bea
&&
\partial_1^{N'_{\lambda,1}}s_\mu([x_1]-[x_2])=
c_\lambda s_{(\mu_1,1^{l'-1})}([x_1]-[x_2]),
\non
\ena
by (iii) of Theorem \ref{th-2-4}. Thus the right hand side of (\ref{eq-3-5})
becomes $c_\lambda \tau_2([x_1]-[x_2])$.
\vskip2mm
\noindent
(iii) This is a direct consequence of Theorem \ref{th-2-4} (iv).
\vskip2mm
\noindent
(iv) Substituting $t=[x_1]-[x_2]$ in $\tau_2(t)$ we have, 
by (iv) of Theorem \ref{th-2-4},
\bea
\tau_2([x_1]-[x_2])&=&
(-1)^{l-1}x_1^{\lambda_1-1}x_2^{l-1}(x_1-x_2)+
\sum \xi_\mu(-1)^{l'-1}x_1^{\mu_1-1}x_2^{l'-1}(x_1-x_2)
\non
\\
&=&
(-1)^{l-1}\bigl(x_1^{\lambda_1}+\sum_\mu \xi_\mu x_1^{\mu_1}\bigr)x_2^{l-1}+O(x_2^l),
\label{eq-3-6}
\ena
where the sum in $\mu$ in the right hand side is over all partitions $\mu$
of the form $\mu=(\mu_1,\lambda_2,...,\lambda_l)$ with $\mu_1>\lambda_1$. Then the term in the bracket
in the right hand side of (\ref{eq-3-6}) is $\tau^{(1)}([x_1])$. 
Thus (iv) is proved. \qed

\section{$\sigma$-function}
In this section we deduce properties of sigma functions from those
of tau functions established in the previous section. To this end we briefly recall the definitions and
properties of sigma functions. 

Let $(n,s)$ be a pair of relatively prime integers satisfying
$2\leq n<s$ and $X$ the compact Riemann surface corresponding to
the algebraic curve defined by
\bea
&&
f(x,y)=0,\qquad
f(x,y)=y^n-x^s-\sum_{ni+sj<ns}\lambda_{ij}x^iy^j.
\label{eq-4-1}
\ena
We assume that the affine curve (\ref{eq-4-1}) is nonsingular.
Then the genus of $X$ is $g=1/2(n-1)(s-1)$. The Riemann surface
$X$ is called an $(n,s)$ curve \cite{BEL2}. 
It has a point $\infty$ over the point $x=\infty$.

For a meromorphic function $F$ on $X$ we denote by $\ord F$ the order
of a pole at $\infty$. The variables $x$ and $y$ can be considered as
meromorphic functions on $X$ which satisfy
\bea
&&
\ord x=n,
\quad
\ord y=s.
\non
\ena
Let $\varphi_i$, $i\geq 1$, be monomials of $x$ and $y$ satisfying
the conditions
\bea
&&
\{\varphi_i\,|\,i\geq 1\}=\{x^iy^j\,|\,i\geq 0,\, n>j\geq 0\},
\non
\\
&&\ord\varphi_i<\ord\varphi_{i+1}, \quad i\geq 1.
\label{eq-4-2}
\ena
For example $\varphi_1=1$, $\varphi_2=x$.

The gap sequence $w_1<\cdots<w_g$ at $\infty$ of $X$ is defined
by
\bea
&&
\{w_i\}={\mathbb Z}_{\geq 0}\backslash \{\ord\varphi_i\,|\,i\geq 1\}.
\non
\ena
It becomes a gap sequence of type $(n,s)$ defined in of Example 1 
in section 2. 

A basis of holomorphic one forms on $X$ is given by
\bea
&&
du_{w_i}:=-\frac{\varphi_{g+1-i}dx}{f_y},
\quad
1\leq i\leq g.
\label{eq-4-3}
\ena
Let $z$ be the local coordinate at $\infty$ such that
\bea
&&
x=\frac{1}{z^n},
\quad
y=\frac{1}{z^s}\left(1+O(z)\right).
\label{eq-4-4}
\ena
Then we have
\bea
&&
du_{w_i}=z^{w_i-1}\left(1+O(z)\right)dz.
\label{eq-4-5}
\ena
We fix an algebraic fundamental form $\hw(p_1,p_2)$ on $X$ \cite{N1}
and decompose it as
\bea
&&
\hw(p_1,p_2)=d_{p_2}\Omega(p_1,p_2)+\sum_{i=1}^g du_{w_i}(p_1)dr_i(p_2),
\non
\ena
where
\bea
&&
\Omega(p_1,p_2)=
\frac{\sum_{i=0}^{n-1}y_1^i[\frac{f(z,w)}{w^{i+1}}]_{+}\vert_{(z,w)=(x_2,y_2)}}
{(x_1-x_2)f_y(x_1,y_1)}dx_1,
\non
\\
&&
[\sum_{m=-\infty}^\infty a_m w^m]_{+}=\sum_{m=0}^\infty a_m w^m.
\non
\ena
Then $dr_i$ automatically becomes a differential of the second kind
whose only singularity is $\infty$ and $\{du_{w_i},dr_i\}$ is a 
symplectic basis of $H^1(X,{\mathbb C})$ \cite{N1}.

We take a symplectic basis of the homology group $H_1(X,{\mathbb Z})$ and define period matrices $\omega_i, \eta_i$, $i=1,2$ by
\bea
&&
2\omega_1=(\int_{\alpha_j}du_{w_i}),
\qquad
2\omega_2=(\int_{\beta_j}du_{w_i}),
\non
\\
&&
-2\eta_1=(\int_{\alpha_j}dr_i),
\qquad
-2\eta_2=(\int_{\beta_j}dr_i).
\non
\ena
The normalized period matrix $\tau$ is given by
$\tau=\omega_1^{-1}\omega_2$.

Let $\tau{}^t\delta'+{}^t\delta''$, $\delta', \delta''\in {\mathbb R}^g$
be the Riemann's constant with respect
to the choice $(\{\alpha_i,\beta_i\},\infty)$, $\delta={}^t(\delta',\delta'')$ and $\theta[\delta](z,\tau)$ the Riemann's theta function
with the characteristic $\delta$.
The sigma function for these data is defined in \cite{BEL1} (see also \cite{N1}).

\begin{defn}\label{def-4-1}
The sigma function $\sigma(u)$, $u={}^t(u_{w_1},...,u_{w_g})$
of an $(n,s)$ curve $X$ is defined by
\bea
&&
\sigma(u)=C\exp\left(\frac{1}{2}{}^tu\eta_1\omega_1^{-1} u\right)
\theta[\delta]((2\omega_1)^{-1}u,\tau),
\non
\ena
where $C$ is a certain constant.
\end{defn}

Let $\lambda=(\lambda_1,...,\lambda_g)$ be the partition corresponding
to the gap sequence at $\infty$ of $X$. Then the constant $C$ is specified
by the condition that the expansion of $\sigma(u)$ at the origin is 
of the form 
\bea
&&
\sigma(u)=s_\lambda(t)\vert_{t_{w_i}=u_{w_i}}+\cdots,
\non
\ena
where $\cdots$ part is a series in $u_{w_i}$ only containing terms 
proportional to $\prod u_{w_i}^{\alpha_i}$ with $\sum \alpha_i w_i>|\lambda|$.

For $m_i\in {\mathbb Z}^g$, $i=1,2$, the sigma function obeys the following
transformation rule:
\bea
&&
\sigma(u+2\sum_{i=1}^2\omega_im_i)
\non
\\
&&
=(-1)^{{}^tm_1m_2+2(\delta'm_1-\delta''m_2)}
\exp\left(2\sum_{i=1}^2 {}^t(\eta_im_i)(u+\sum_{i=1}^2 \omega_im_i)\right)
\sigma(u).
\label{eq-4-6}
\ena

Let $A$ be the affine ring of $X\backslash\{\infty\}$. As a vector space
$\{\varphi_i|i\geq 1\}$ is a basis of $A$. We embed $A$ into UGM using the local coordinate $z$ as in
\cite{N2}. Then the tau function $\tau(t)$ of the KP-hierarchy corresponding to this point of UGM has the form
\bea
&&
\tau(t)=s_\lambda(t)+\sum_{\lambda<\mu}\xi_\mu s_\mu(t).
\non
\ena
It can be expressed in terms of the sigma function as
\bea
&&
\tau(t)=
\exp\left(-\sum_{i=1}^\infty c_it_i+\frac{1}{2}{\widehat q}(t)\right)
\sigma(Bt),
\label{eq-4-7}
\ena
where ${\widehat q}(t)=\sum {\widehat q}_{ij}t_it_j$, 
$B=(b_{ij})_{1\leq i\leq g, 1\leq j}$ a certain $g\times \infty$ 
matrix satisfying the condition
\bea
&&
b_{ij}=\left\{
\begin{array}{cc}
0& \quad\hbox{if $j<w_i$}\\
1& \quad\hbox{if $j=w_i$},
\end{array}
\right.
\label{eq-4-8}
\ena
and $c_i, {\widehat q}_{ij}, b_{ij}$ are certain constants \cite{N2}
\footnote{In the defining equation of $c_i$ in \cite{N2} $c_i$ should be corrected to $c_i/i$} \cite{EEG}\cite{EH}. The constant $c_i$ are irrelevant to $c_k$ in 
Theorem \ref{th-2-1} and is not used in other parts of this paper.

In this section $\partial_i$ is used for 
$\partial/\partial t_i$ as in the previous 
section and $\partial_{u_i}$ is used for 
$\partial/\partial u_i$.

A point $p\in X$ is identified with its Abel-Jacobi image 
$\displaystyle{\int_\infty^pd{\mathbf u}}$, where $d{\mathbf u}
={}^t(du_{w_1},...,du_{w_g})$. 
By the definition of the matrix $B$, for $p\in X$, 
the following equation is valid:
\bea
&&
B[z(p)]=p,
\label{eq-4-9}
\ena
where $z(p)$ is the value of the local coordinate $z$ at $p$ and $[z(p)]=[z(p),z(p)^2/2,...]$ as before.

Corresponding to Proposition \ref{prop-3-1} we have

\begin{theorem}\label{th-4-1}
Let $0\leq k\leq g-1$. If $\sum_{i=1}^g\alpha_i w_i<N_{\lambda,k}$ them
\bea
&&
\partial_{u_{w_1}}^{\alpha_1}\cdots \partial_{u_{w_g}}^{\alpha_g}
\sigma\bigl(\,\sum_{i=1}^k p_i\bigr)=0,
\non
\ena
for $p_1,...,p_k\in X$.
\end{theorem}
\vskip3mm

\noindent
{\bf Remark 2} In the case of the curve $y^n=f(x)$ Theorem \ref{th-4-1}
is proved in \cite{MP}.
\vskip2mm
\begin{lemma}\label{lem-4-1}
Let $0\leq k\leq g-1$. If $\wt\alpha<N_{\lambda,k}$ 
\bea
&&
\partial^\alpha\sigma(Bt)\vert_{t=t^{(k)}}=0,
\non
\ena
where $t^{(k)}=\sum_{i=1}^k[z_i]$, 
$z_i=z(p_i)$ and $p_1,...,p_k\in X$.
\end{lemma}
\vskip2mm
\noindent
\pf The assertion easily follows from (\ref{eq-4-7}) and
Proposition \ref{prop-3-1}. \qed
\vskip2mm

\noindent
{\it Proof of Theorem \ref{th-4-1}.}\par
We introduce the lexicographical order on ${\mathbb Z}_{\geq 0}^g$
comparing from the right. Namely define $(\alpha_1,...,\alpha_g)<(\beta_1,...,\beta_g)$ if there exists $1\leq i\leq g$ such that 
$\alpha_g=\beta_g$,...,$\alpha_{i+1}=\beta_{i+1}$ and 
$\alpha_i<\beta_i$.

We prove
\bea
&&
\partial_{u_{w_1}}^{\beta_1}\cdots
\partial_{u_{w_g}}^{\beta_g}
\sigma(Bt^{(k)})=0,
\qquad
\sum_{i=1}^g \beta_iw_i<N_{\lambda,k},
\label{eq-4-11}
\ena
by induction on the order of $(\beta_1,...,\beta_g)$.

The case $(\beta_1,...,\beta_g)=(0,...,0)$ is obvious by Lemma \ref{lem-4-1}.

Take $(\beta_1,...,\beta_g)>(0,...,0)$ such that $\sum_{i=1}^g \beta_iw_i<N_{\lambda,k}$. Suppose that (\ref{eq-4-11})
is valid for any $(\beta_1',...,\beta_g')$ satisfying 
$(\beta_1',...,\beta_g')<(\beta_1,...,\beta_g)$ and $\sum_{i=1}^g \beta_i'w_i<N_{\lambda,k}$.

Notice that $\sigma(Bt)$ is a composition of $\sigma(u)$ with
\bea
&&
u_{w_i}=t_{w_i}+\sum_{w_i<j}b_{ij} t_j,
\quad
1\leq i\leq g.
\label{eq-4-10}
\ena
By the chain rule, 
\bea
&&
\partial_{w_i}=\partial_{u_{w_i}}+\sum_{l<i}b_{lw_i}\partial_{u_{w_l}}.
\label{eq-4-12}
\ena

Let $j$ be the maximum number such that $\beta_j\neq 0$.
Then 
\bea
&&
\partial_{w_1}^{\beta_1}\cdots\partial_{w_{g}}^{\beta_{g}} \sigma(Bt)
\non
\\
&=&
\partial_{u_1}^{\beta_1}(\partial_{u_{w_2}}+b_{1w_2}\partial_{u_1})^{\beta_2}
\cdots(\partial_{u_{w_j}}+\sum_{l<j}b_{lw_j}\partial_{u_{w_l}})^{\beta_j}
\sigma(Bt)
\non
\\
&=&
\partial_{u_1}^{\beta_1}\partial_{u_{w_2}}^{\beta_2}
\cdots\partial_{u_{w_j}}^{\beta_j}
\sigma(Bt)+\cdots,
\label{eq-4-13}
\ena
where $\cdots$ part contains terms of the form
\bea
&&
\partial_{u_1}^{\gamma_1}\cdots\partial_{u_{w_g}}^{\gamma_g}\sigma(Bt),
\quad
\sum_{i=1}^g \gamma_iw_i<N_{\lambda,k},
\quad
(\gamma_1,...,\gamma_g)<(\beta_1,...,\beta_g).
\non
\ena
At $t=t^{(k)}$ the left hand side of (\ref{eq-4-13}) vanishes by 
Lemma \ref{lem-4-1} and $\cdots$ part in the right hand side of (\ref{eq-4-13})
vanishes by the assumption of induction.
Thus (\ref{eq-4-11}) is proved. \qed

Corresponding to Theorem \ref{th-3-1} and Corollary \ref{cor-3-1} we have

\begin{cor}\label{cor-4-1}
Let $1\leq k\leq g$  and $p_1,...,p_k\in X$. Then
\vskip2mm
\noindent
(i) We have  
\bea
&&
\partial_{u_{a^{(k)}_1}}\cdots\partial_{u_{a^{(k)}_{m_k}}}
\sigma\bigl(\,\sum_{i=1}^k p_i\,\bigr)
=c_k S_{(\lambda_1,...,\lambda_k)}
(z_1,...,z_k)+\cdots,
\non
\ena
where $\cdots$ part is a series in $z_i$ containing only terms
proportional to $\prod_{i=1}^k z_i^{\alpha_i}$ with 
$\sum_{i=1}^k \alpha_i>\sum_{i=1}^k\lambda_i$.
\vskip2mm

\noindent
(ii) The following expansion is valid:
\bea
&&
\partial_{u_{a^{(k)}_1}}\cdots\partial_{u_{a^{(k)}_{m_k}}}
\sigma\bigl(\,\sum_{i=1}^k p_i\,\bigr)
=\frac{c_k}{c_{k-1}}
\partial_{u_{a^{(k-1)}_1}}\cdots\partial_{u_{a^{(k-1)}_{m_{k-1}}}}
\sigma\bigl(\,\sum_{i=1}^{k-1} p_i\,\bigr)
z_k^{\lambda_k}+O(z_k^{\lambda_{k}+1}).
\non
\ena
\end{cor}
\vskip2mm
\noindent
\pf By Theorem \ref{th-4-1} and (\ref{eq-4-12}) we have
\bea
&&
\partial_{a^{(k)}_1}\cdots\partial_{a^{(k)}_{m_k}}\sigma(Bt)\vert_{t=t^{(k)}}
=
\partial_{u_{a^{(k)}_1}}\cdots\partial_{u_{a^{(k)}_{m_k}}}
\sigma\bigl(\,\sum_{i=1}^k p_i\,\bigr).
\label{eq-4-14}
\ena
Let us write (\ref{eq-4-7}) as $\sigma(Bt)=\varepsilon(t)\tau(t)$ with
\bea
&&
\varepsilon(t)=
\exp\left(\sum_{i=1}^\infty c_it_i-\frac{1}{2}{\widehat q}(t)\right).
\non
\ena
By Proposition \ref{prop-3-1} and Corollary \ref{cor-3-1} we have
\bea
&&
\partial_{a^{(k)}_1}\cdots\partial_{a^{(k)}_{m_k}}\sigma(Bt)\vert_{t=t^{(k)}}
\non
\\
&=&
\varepsilon(t^{(k)})
\partial_{a^{(k)}_1}\cdots\partial_{a^{(k)}_{m_k}}
\tau(t^{(k)})
\non
\\
&=&
c_{k-1}^{-1}c_k 
\varepsilon(t^{(k)})
\partial_{a^{(k-1)}_1}\cdots\partial_{a^{(k-1)}_{m_{k-1}}}
\tau(t^{(k-1)})z_k^{\lambda_k}+O(z_k^{\lambda_k+1})
\non
\\
&=&
c_{k-1}^{-1}c_k 
\partial_{a^{(k-1)}_1}\cdots\partial_{a^{(k-1)}_{m_{k-1}}}
\sigma(Bt)\vert_{t=t^{(k-1)}}z_k^{\lambda_k}+O(z_k^{\lambda_k+1}).
\label{eq-4-14'}
\ena
Then the assertion (ii) follows from (\ref{eq-4-14}) and
the assertion (i) follows from the second line of (\ref{eq-4-14'})
, Theorem \ref{th-3-1} (i) and the definition of $\tau^{(k)}(t)$. 
\qed

The following corollary can similarly be proved using  
Theorem \ref{th-3-2} and Corollary \ref{cor-3-2}.

\begin{cor}\label{cor-4-2}
Let $1\leq k\leq g$ and $p_1,...,p_k\in X$. Then
\vskip2mm
\noindent
(i) We have 
\bea
&&
\partial_{u_1}^{N_{\lambda,k}}\sigma\bigl(\,\sum_{i=1}^k p_i\,\bigr)=
c'_{\lambda,k}S_{(\lambda_1,...,\lambda_k)}(z_1,...,z_k)+\cdots,
\non
\ena
where $\cdots$ part is a series in $z_i$ containing only terms
proportional to $\prod_{i=1}^k z_i^{\alpha_i}$ with 
$\sum_{i=1}^k \alpha_i>\sum_{i=1}^k\lambda_i$.
\vskip2mm
\noindent
(ii) The following expansion holds:
\bea
&&
\partial_{u_1}^{N_{\lambda,k}}
\sigma\bigl(\,\sum_{i=1}^k p_i\,\bigr)
=
\frac{c'_{\lambda,k}}{c'_{\lambda,k-1}}
\partial_{u_1}^{N_{\lambda,k-1}}\sigma
\bigl(\,\sum_{i=1}^{k-1} p_i\,\bigr)
z_k^{\lambda_k}
+O(z_k^{\lambda_k+1}).
\non
\ena
\end{cor}
\vskip2mm
\noindent

Corresponding to Theorem \ref{th-3-3} we have

\begin{theorem}\label{th-4-2}
(i) If $n<N'_{\lambda,1}$ we have, for $p_1,p_2\in X$, 
\bea
&&
\partial_{u_1}^n\sigma(p_1-p_2)=0.
\non
\ena
\vskip2mm
\noindent
(ii) The following expansion with respect to $z_i=z(p_i)$, $i=1,2$
is valid:
\bea
&&
\partial_{u_1}^{N'_{\lambda,1}} \sigma(p_1-p_2)
=(-1)^{g-1} c_{\lambda} (z_1z_2)^{g-1}(z_1-z_2)(1+\cdots),
\non
\ena
where $\cdots$ part is a series in $z_1,z_2$ which 
contains only terms proportional to $z_1^iz_2^j$ with $i+j>0$.
\vskip2mm
\noindent
(iii) We have
\bea
&&
\partial_{u_1}^{N'_{\lambda,1}}\sigma(p_1-p_2)=
(-1)^{g-1}
\frac{c_\lambda}{c'_{\lambda,1}}
\partial_{u_1}^{N_{\lambda,1}}\sigma(p_1)z_2^{g-1}+O(z_2^g).
\non
\ena
\end{theorem}
\vskip2mm
\noindent
\pf (i) Notice that 
\bea
&&
\partial_1^m\sigma(Bt)=\partial_{u_1}^m\sigma(Bt)
\label{eq-4-15}
\ena
for any $m$. 
Differentiating $\sigma(Bt)=\varepsilon(t)\tau(t)$ and using (\ref{eq-4-15})
and (i) of Theorem \ref{th-3-3} we have the assertion.
\vskip2mm
\noindent
(ii) We have, by (i), (ii), (iii) of Theorem \ref{th-3-3},
\bea
\partial_{u_1}^{N'_{\lambda,1}}\sigma(p_1-p_2)&=&
\partial_1^{N'_{\lambda,1}}\sigma(Bt)\vert_{t=[z_1]-[z_2]}
\non
\\
&=&
\varepsilon([z_1]-[z_2])\partial_1^{N'_{\lambda,1}}
\tau([z_1]-[z_2])
\non
\\
&=&
c_\lambda\varepsilon([z_1]-[z_2])\tau_2([z_1]-[z_2])
\non
\\
&=&c_\lambda(-1)^{g-1}(z_1z_2)^{g-1}(z_1-z_2)(1+\cdots).
\non
\ena

\noindent
(iii) In the computation in (ii) we have
\bea
c_\lambda\varepsilon([z_1]-[z_2])\tau_2([z_1]-[z_2])
&=&(-1)^{g-1}c_\lambda (c'_{\lambda,1})^{-1}
\varepsilon([z_1]-[z_2])
\partial_1^{N_{\lambda,1}}\tau([z_1])z_2^{g-1}+O(z_2^g)
\non
\\
&=&
(-1)^{g-1}c_\lambda (c'_{\lambda,1})^{-1}
\partial_1^{N_{\lambda,1}}
\sigma(Bt)\vert_{t=[z_1]} z_2^{g-1}+O(z_2^g),
\non
\ena
by Theorem \ref{th-3-3} (i), (ii), (iv) and Theorem \ref{th-3-2}.
 Then the assertion follows from (\ref{eq-4-15}). \qed

\section{Addition formulae}
Let $E(p_1,p_2)$ be the prime form \cite{F,M2} of an $(n,s)$ curve $X$. 
In \cite{N1} we have introduced the prime function ${\tilde E}(p_1,p_2)$ by
\bea
&&
{\tilde E}(p_1,p_2)=-E(p_1,p_2)\prod_{i=1}^2\sqrt{du_{w_g}(p_i)}
\exp\left(\frac{1}{2}\int_{p_1}^{p_2}{}^td{\mathbf u}\eta_1\omega_1^{-1}
\int_{p_1}^{p_2}d{\mathbf u}\right).
\label{eq-5-1}
\ena
Notice that ${\tilde E}(p_1,p_2)$ is not a $(-1/2,-1/2)$ form but 
a multi-valued analytic function on $X$
and thus it has a sense to talk about the transformation rule 
if $p_i$ goes around a cycle of $X$.

The prime function has the following properties.
\vskip2mm
\noindent
(i) ${\tilde E}(p_2,p_1)=-{\tilde E}(p_1,p_2)$.
\vskip2mm
\noindent
(ii) As a function of $p_1$, the zero divisor of
${\tilde E}(p_1,p_2)$ is $p_2+(g-1)\infty$.
\vskip2mm
\noindent
(iii) Let $m_i={}^t(m_{i1},...,m_{ig})\in {\mathbb Z}^g$.
If $p_2$ goes round the cycle 
$\gamma=\sum_{i=1}^g(m_{1i}\alpha_i+m_{2i}\beta_i)$, ${\tilde E}(p_1,p_2)$
transforms as
\bea
&&
{\tilde E}(p_1,\gamma p_2)=T(m_1,m_2\vert \int_{p_1}^{p_2}d{\mathbf u})
\, {\tilde E}(p_1,p_2),
\label{eq-5-2}
\ena
with
\bea
T(m_1,m_2\vert u)=
(-1)^{{}^tm_1m_2+2(\delta'm_1-\delta''m_2)}
\exp\bigl(2\sum_{i=1}^2{}^t(\eta_im_i)(u+\sum_{i=1}^2\omega_im_i)\bigr).
\non
\ena

\noindent
(iv) At $(\infty,\infty)$, ${\tilde E}(p_1,p_2)$ has the expansion
of the form
\bea
&&
{\tilde E}(p_1,p_2)=(z_1z_2)^{g-1}(z_1-z_2)(1+\sum_{i+j\geq 1}c_{ij}
z_1^iz_2^j),
\label{eq-5-3}
\ena
where $z_i=z(p_i)$.
\vskip2mm

The specialization ${\tilde E}(\infty,p)$ of ${\tilde E}(p_1,p_2)$ 
is defined by 
\bea
&&
-{\tilde E}(p_1,p_2)={\tilde E}(\infty,p_2)z_1^{g-1}+O(z_1^g).
\label{eq-5-4}
\ena
It has the following properties corresponding to (iii) and (iv) above.
\vskip2mm

\noindent
(iii)' Under the same notation as in (iii) for ${\tilde E}(p_1,p_2)$ we have
\bea
&&
{\tilde E}(\infty,\gamma p_2)=
T(m_1,m_2\vert \int_{\infty}^{p_2}d{\mathbf u})\,{\tilde E}(\infty,p_2).
\label{eq-5-5}
\ena

\noindent
(iv)' ${\tilde E}(\infty, p_2)=z_2^g+O(z_2^{g+1})$.
\vskip2mm

The following theorem gives the expression of the prime function
in terms of a derivative of the sigma function.

\begin{theorem}\label{th-5-1}
Let $\lambda=(\lambda_1,...,\lambda_g)$ be the partition corresponding
to the gap sequence at $\infty$ of an $(n,s)$ curve $X$ .
Then 
\bea
&&
{\tilde E}(p_1,p_2)=(-1)^{g-1}c_\lambda^{-1}\partial_{u_1}^{N'_{\lambda,1}}
\sigma(p_1-p_2).
\non
\ena
\end{theorem}

\begin{lemma}\label{lem-5-1}
Under the same notation as in (\ref{eq-5-2}) we have
\bea
&&
\partial_{u_1}^{N'_{\lambda,1}}
\sigma(p_1-\gamma p_2)
=T(m_1,m_2\vert \int_{p_1}^{p_2}d{\mathbf u})\,
\partial_{u_1}^{N'_{\lambda,1}}
\sigma(p_1-p_2).
\non
\ena
\end{lemma}
\vskip2mm
\noindent
\pf 
Notice that
$\gamma p_2=p_2+2\omega_1m_1+2\omega_2m_2$ and
\bea
&&
\sigma(u-2\sum_{i=1}^2 \omega_im_i)=T(-m_1,-m_2|u)\sigma(u).
\label{eq-5-6}
\ena
Applying $\partial_{u_1}^{N'_{\lambda,1}}$ to (\ref{eq-5-6})
and setting $u=p_1-p_2$, we get, by Theorem \ref{th-4-2} (i)
we have
\bea
&&
\partial_{u_1}^{N'_{\lambda,1}}\sigma(p_1-p_2-2\sum_{i=1}^2 \omega_im_i)=
T(-m_1,-m_2|\int_{p_2}^{p_1}d{\mathbf u})\,
\partial_{u_1}^{N'_{\lambda,1}}\sigma(p_1-p_2).
\label{eq-5-7}
\ena
Then the assertion follows from $T(-m_1,-m_2|u)=T(m_1,m_2|-u)$. \qed

\vskip5mm
\noindent
{\it Proof of Theorem \ref{th-5-1}.}\par
Notice that 
\bea
&&
\partial_{u_1}^{N'_{\lambda,1}}\sigma(-u)=-
\partial_{u_1}^{N'_{\lambda,1}}\sigma(u),
\label{eq-5-8}
\ena
since $\sigma(-u)=(-1)^{|\lambda|}\sigma(u)$ \cite{N1} and 
$N'_{\lambda,1}=|\lambda|-2g+1$.

Consider the function
\bea
&&
F(p_1,p_2)=\frac{\partial_{u_1}^{N'_{\lambda,1}}
\sigma(p_1-p_2)}{{\tilde E}(p_1,p_2)}.
\label{eq-5-9}
\ena
It is symmetric in $p_1$ and $p_2$ by (\ref{eq-5-8}), (i) of properties of 
${\tilde E}(p_1,p_2)$ and is a meromorphic function on $X\times X$ by Lemma \ref{lem-5-1}. Fix $p_1$ near $\infty$. As a function of $p_2$ $F(p_1,p_2)$ has no
singularity by Theorem \ref{th-4-1}, Theorem \ref{th-4-2} (ii)  and the property (ii) of ${\tilde E}(p_1,p_2)$. 
Therefore it does not depend on $p_2$. 
It means that, for some 
non-empty open neighborhood $U$ of $\infty$, $F(p_1,p_2)$ does not depend on $p_2$ on $U\times X$. 
Since $F(p_1,p_2)$ is symmetric, it is a constant on 
$U\times U$. 
Thus it is a constant on $X\times X$ because it is meromorphic. 
The constant can be determined by comparing the expansion using Theorem \ref{th-4-2} (ii) and the property (iv) of ${\tilde E}(p_1,p_2)$. \qed

\begin{cor}\label{cor-5-1}
For $p\in X$ we have
\bea
&&
{\tilde E}(\infty,p)=
{c'_{\lambda,1}}^{-1}\partial_{u_1}^{N_{\lambda,1}}\sigma(p).
\non
\ena
\end{cor}
\vskip2mm
\noindent
\pf Compare the expansion in the equation of Theorem \ref{th-5-1} using (iii) of Theorem \ref{th-4-2}. \qed
\vskip2mm

\vskip2mm
\noindent
{\bf Remark 3} In the case of hyperelliptic curves the prime function
can be given using the derivative determined from the sequence $a^{(2)}_j$.
This is because $p_1-p_2$ can be written as a sum $p_1+p_2^\ast$ where ${}^\ast$denoting the hyperelliptic involution. Such expression is given in \cite{GMO}.
\vskip2mm

The following theorem had been proved in \cite{N1}.

\begin{theorem}\cite{N1}\label{th-5-2}
For $n\geq g$ and $p_i\in X$, $1\leq i\leq n$,
\bea
&&
\sigma\bigl(\,\sum_{i=1}^n p_i\,\bigr)=
\frac{\prod_{i=1}^n {\tilde E}(\infty,p_i)^n}{\prod_{i<j}^n{\tilde E}(p_i,p_j)}
\det\left(\varphi_i(p_j)\right)_{1\leq i,j\leq n}.
\non
\ena
\end{theorem}

By comparing the top term of the series expansion in $z(p_n)$,
using Theorem \ref{th-2-2} and Corollary \ref{cor-4-2}, beginning from
$n=g$ successively in the equation of this theorem we get

\begin{cor}\label{cor-5-2}
For $n<g$ we have
\bea
&&
\partial_{u_1}^{N_{\lambda,n}}
\sigma\bigl(\,\sum_{i=1}^n p_i\,\bigr)=
c'_{\lambda,n}
\frac{\prod_{i=1}^n {\tilde E}(\infty,p_i)^n}{\prod_{i<j}^n{\tilde E}(p_i,p_j)}
\det\left(\varphi_i(p_j)\right)_{1\leq i,j\leq n}.
\non
\ena
\end{cor}

Combining Theorem \ref{th-5-1} and Corollaries \ref{cor-5-1} and \ref{cor-5-2}
we have the following addition formulae for sigma functions.

\begin{cor}\label{cor-5-3}
(i) For $n\geq g$ and $p_i\in X$, $1\leq i\leq n$,
\bea
&&
\frac{\sigma\bigl(\,\sum_{i=1}^n p_i\,\bigr)
\prod_{i<j}\partial_{u_1}^{N'_{\lambda,1}}
\sigma(p_j-p_i)}{\prod_{i=1}^n(\partial_{u_1}^{N_{\lambda,1}}\sigma(p_i))^n}
=b_{\lambda,n}
\det\left(\varphi_{i}(p_j)\right)_{1\leq i,j\leq n},
\non
\ena
with
\bea
&&
b_{\lambda,n}=
(-1)^{\frac{1}{2}gn(n-1)}c_{\lambda}^{\frac{1}{2}n(n-1)}(c'_{\lambda,1})^{-n^2}.
\non
\ena
\noindent
(ii) For $n<g$ 
\bea
&&
\frac{\partial_{u_1}^{N_{\lambda,n}}
\sigma\bigl(\,\sum_{i=1}^n p_i\,\bigr)
\prod_{i<j}\partial_{u_1}^{N'_{\lambda,1}}\sigma(p_j-p_i)}
{\prod_{i=1}^n(\partial_{u_1}^{N_{\lambda,1}}\sigma(p_i))^n}
=b_{\lambda,n}
\det\left(\varphi_{i}(p_j)\right)_{1\leq i,j\leq n},
\non
\ena
with
\bea
&&
b_{\lambda,n}=
(-1)^{\frac{1}{2}gn(n-1)}c_{\lambda}^{\frac{1}{2}n(n-1)}(c'_{\lambda,1})^{-n^2}
c'_{\lambda,n}.
\non
\ena
\end{cor}

Similarly, using Theorem \ref{th-5-2}, Theorem \ref{th-2-1} and Corollary \ref{cor-4-1}, we have

\begin{cor}\label{cor-5-4}
For $n<g$ and $p_i\in X$, $1\leq i\leq n$, we have
\bea
&&
\partial_{u_{a^{(n)}_1}}\cdots\partial_{u_{a^{(n)}_{m_n}}}
\sigma(\sum_{i=1}^n p_i)
=c_n
\frac{\prod_{i=1}^n {\tilde E}(\infty,p_i)^n}{\prod_{i<j}^n{\tilde E}(p_i,p_j)}
\det\left(\varphi_i(p_j)\right)_{1\leq i,j\leq n}.
\non
\ena
\end{cor}

\begin{cor}\label{cor-5-5}
For $n<g$ and $p_i\in X$, $1\leq i\leq n$, we have
\bea
&&
\frac{\partial_{u_{a^{(n)}_1}}\cdots\partial_{u_{a^{(n)}_{m_n}}}
\sigma(\sum_{i=1}^n p_i)\prod_{i<j}\partial_{u_1}^{N'_{\lambda,1}}
\sigma(p_j-p_i)}{\prod_{i=1}^n(\partial_{u_1}^{N_{\lambda,1}}\sigma(p_i))^n}
=b'_{\lambda,n}
\det\left(\varphi_{i}(p_j)\right)_{1\leq i,j\leq n},
\non
\ena
with
\bea
&&
b'_{\lambda,n}=
(-1)^{\frac{1}{2}gn(n-1)}c_{\lambda}^{\frac{1}{2}n(n-1)}
(c'_{\lambda,1})^{-n^2}c_n.
\non
\ena
\end{cor}

In the case of hyperelliptic curves $\partial_{u_1}^{N'_{\lambda,1}}
\sigma(p_j-p_i)$ can be replaced by a constant multiple of 
"$a^{(2)}_j$-derivative" as remarked before (Remark after Corollary \ref{cor-5-1}).
Then  Corollaries \ref{cor-5-3}, \ref{cor-5-5} 
recovers the formulae in \cite{O1}.
\vskip10mm

\noindent
{\large {\bf Acknowledgments}} 
\vskip3mm
\noindent
The authors would like to thank Yasuhiko
Yamada for suggesting Theorem \ref{th-2-3} and Yoshihiro \^Onishi
for a stimulating discussion. 
The first author is grateful to the organizers of the conference 
"Symmetries, Integrable Systems and Representations" held at
Lyon, December 13-16, 2011, for financial support and kind hospitality.
This research is partially supported by JSPS Grant-in-Aid for Scientific Research (C) No.23540245.
\vskip3mm

\noindent
{{\bf Note}}: A part of the results in the present paper is reported
 in \cite{Y}.

\end{document}